\documentclass [12pt,a4paper] {amsart}

\usepackage [T1] {fontenc}
\usepackage {amssymb}
\usepackage {amsmath,bm}
\usepackage {amsthm}
\usepackage{tikz}
\usepackage{hyperref}
\usepackage{url}
\usepackage{soul}

 \usepackage[titletoc,toc,title]{appendix}

\usetikzlibrary{arrows,decorations.pathmorphing,backgrounds,positioning,fit,petri}
\usetikzlibrary{matrix}

\usepackage{enumitem}

\def\MLine#1{\par\hspace*{-\leftmargin}\parbox{\textwidth}{\[#1\]}}

\DeclareMathOperator {\td} {td}

\DeclareMathOperator {\Th} {Th}

\DeclareMathOperator {\tp} {tp}

\DeclareMathOperator {\ord} {ord}
\DeclareMathOperator {\acl} {acl}

\DeclareMathOperator {\cl} {cl}
\DeclareMathOperator {\im} {Im}
\DeclareMathOperator {\pr} {pr}

\DeclareMathOperator {\DCF} {DCF}

\DeclareMathOperator {\Mr} {MR}
\DeclareMathOperator {\MD} {MD}
\DeclareMathOperator {\SL} {SL}

\DeclareMathOperator {\GL} {GL}

\DeclareMathOperator {\Loc} {Loc}

\DeclareMathOperator {\alg} {alg}

\DeclareMathOperator {\ppr} {Pr}

\DeclareMathOperator {\A} {\mathsf{A}}
\DeclareMathOperator {\B} {\mathsf{B}}
\DeclareMathOperator {\K} {\mathsf{K}}

\DeclareMathOperator {\X} {\mathsf{X}}
\DeclareMathOperator {\Y} {\mathsf{Y}}
\DeclareMathOperator {\C} {\mathsf{C}}

\newcommand {\liff} {\leftrightarrow}
\newcommand {\seq} {\subseteq}

\theoremstyle {definition}
\newtheorem {definition}{Definition} [section]

\theoremstyle {plain}

\newtheorem {lemma} [definition] {Lemma}
\newtheorem {theorem} [definition] {Theorem}
\newtheorem {proposition} [definition] {Proposition}
\newtheorem {corollary} [definition] {Corollary}
\newtheorem {conjecture} [definition] {Conjecture}

\newtheorem{innercustomthm}{Theorem}
\newenvironment{customthm}[1]
  {\renewcommand\theinnercustomthm{#1}\innercustomthm}
  {\endinnercustomthm}

\theoremstyle {remark}
\newtheorem {remark} [definition] {Remark}

\makeatletter
\newcommand {\forksym} {\raise0.2ex\hbox{\ooalign{\hidewidth$\vert$\hidewidth\cr\raise-0.9ex\hbox{$\smile$}}}}
\def\@forksym@#1#2{\mathrel{\mathop{\forksym}\displaylimits_{#2}}}
\def\forkind{\@ifnextchar_{\@forksym@}{\forksym}}

\makeatother

\calclayout

\begin {document}

\title{Ax-Schanuel and strong minimality for the $j$-function}

\author{Vahagn Aslanyan}
\address{Department of Mathematical Sciences, CMU, 5000 Forbes Ave, Pittsburgh, PA 15213, USA}
\address{Institute of Mathematics of the National Academy of Sciences of Armenia, Yerevan 0019, Armenia}
\address{\texttt{Current address:} School of Mathematics, University of East Anglia, Norwich, NR4 7TJ, UK}
\email{V.Aslanyan@uea.ac.uk}

\thanks{This work was done while the author was a DPhil student at the University of Oxford, and later a postdoctoral associate at Carnegie Mellon University. Some revisions were made at the University of East Anglia. Partially supported by Dulverton Scholarship, Luys Scholarship and AGBU UK Scholarship (at the University of Oxford) and EPSRC grant EP/S017313/1 (at the University of East Anglia).}



\maketitle

\begin{abstract}
Let $\mathcal{K}:=(K;+,\cdot, D, 0, 1)$ be a differentially closed field of characteristic $0$ with field of constants $C$.

In the first part of the paper we explore the connection between Ax-Schanuel type theorems (predimension inequalities) for a differential equation $E(x,y)$ and the geometry of the fibres $U_s:=\{ y:E(s,y) \wedge y \notin C \}$ where $s$ is a non-constant element. We show that certain types of predimension inequalities imply strong minimality and geometric triviality of $U_s$. Moreover, the induced structure on the Cartesian powers of $U_s$ is given by special subvarieties. 
In particular, since the $j$-function satisfies an Ax-Schanuel inequality of the required form (due to Pila and Tsimerman),  {applying our results to the $j$-function we recover} a theorem of Freitag and Scanlon stating that the differential equation of $j$ defines a strongly minimal set with trivial geometry.

In the second part of the paper we study strongly minimal sets in the $j$-reducts of differentially closed fields. Let $E_j(x,y)$ be the (two-variable) differential equation of the $j$-function. We prove a Zilber style classification result for strongly minimal sets in the reduct $\K:=(K;+, \cdot, E_j)$. More precisely, we show that in $\K$ all strongly minimal sets are geometrically trivial or non-orthogonal to $C$.  Our proof is based on the Ax-Schanuel theorem and a matching \emph{Existential Closedness} statement which asserts that systems of equations in terms of $E_j$ have solutions in $\K$ unless having a solution contradicts Ax-Schanuel.
\end{abstract}

\section{Introduction}

Throughout the paper we let $\mathcal{K}=(K;+,\cdot, D, 0, 1)$ be a differentially closed field of characteristic $0$ with field of constants $C$. 

Let $E(x,y)$ be the set of solutions of a differential equation $f(x,y)=0$ with rational or, more generally, constant coefficients. A general question that we are interested in is whether $E$ satisfies an Ax-Schanuel type inequality. A motivating example is the exponential differential equation $Dy=yDx$. We know that the original Ax-Schanuel theorem \cite{Ax} gives a predimension inequality (in the sense of Hrushovski \cite{Hru}) which governs the geometry of our equation. In this case the corresponding reduct of a differentially closed field can be reconstructed by a Hrushovski-style amalgamation-with-predimension construction \cite{Kirby-semiab}. This kind of predimension inequalities is said to be \emph{adequate} (see \cite{Aslanyan-thesis,Aslanyan-adequate-predim} for a precise definition). This means that the reduct satisfies an \emph{existential closedness} property which asserts roughly that a system of exponential differential equations which is not overdetermined has a solution. Being overdetermined means that the existence of a solution would contradict Ax-Schanuel. Thus, having an adequate Ax-Schanuel type inequality for $E$ would give us a good understanding of its model theory. For more details on this and related problems see \cite{Aslanyan-thesis,Aslanyan-adequate-predim,Kirby-semiab,Zilb-pseudoexp,Zilb-analytic-pseudoanalytic}. 
Moreover, Ax-Schanuel type statements have applications in diophantine geometry. Indeed, Ax-Schanuel and its weak versions play an important role in the proofs of many diophantine theorems related to unlikely intersections and the famous Zilber-Pink conjecture (see for example \cite{Zilb-exp-sum-published,Bom-Mas-Zan,Pila-Tsim-Ax-j,Kirby-semiab,Habegger-Pila-o-min-certain,Daw-Ren,Aslanyan-weakMZPD}).

Thus, we want to classify differential equations in two variables with respect to the property of satisfying an Ax-Schanuel type inequality. The present work should be seen as a part of that more general project. In the first part of the paper we explore the connection between Ax-Schanuel type theorems (predimension inequalities) for a differential equation $E(x,y)$ and the geometry of the fibres of $E$. More precisely, given a predimension inequality (not necessarily adequate) for solutions of $E$ of a certain type, we show that the fibres of $E$ are strongly minimal and geometrically trivial after removing constant points. Moreover, the induced structure on the Cartesian powers of those fibres is given by \emph{special} subvarieties.

One of the main results of the first part is as follows (for the definition of weakly $\mathcal{P}$-special varieties see Section \ref{j-setup}).

\begin{theorem}\label{thm:intro-AS-SM}
Let $E(x,y)$ be defined by $R(x, y, \partial_x y, \ldots, \partial_x^m y)=0$ where $\partial_x = \frac{D}{Dx}$ for a non-constant $x$ and $R(X,\bar{Y})$ is an algebraic polynomial over $C$, irreducible over $C(X)^{\alg}$. Assume $E$ satisfies the following Ax-Schanuel condition for a collection $\mathcal{P}$ of algebraic polynomials $P(X,Y) \in C[X,Y]$:
\begin{itemize}[leftmargin=0.4cm]
\item[] Let $x_1,\ldots,x_n,y_1,\ldots,y_n$ be non-constant elements of $K$ with $E(x_i,y_i)$. If $P(y_i,y_j)\neq 0$ for all $P\in\mathcal{P}$ and $i \neq j$ then 
\MLine {\td_CC(x_1, y_1, \partial_{x_1}y_1, \ldots, \partial^{m-1}_{x_1}y_1, \ldots, x_n, y_n, \partial_{x_n}y_n, \ldots, \partial^{m-1}_{x_n}y_n) \geq mn+1.}
\end{itemize}
Then for every $s\in K\setminus C$ the set $U_s:=\{y: E(s,y) \wedge y\notin C \}$ is strongly minimal with trivial geometry. Furthermore, every definable subset of a Cartesian power of $U_s$ is a Boolean combination of weakly $\mathcal{P}$-special subvarieties.
\end{theorem}

In particular, let $F_j(y,Dy,D^2y,D^3y)=0$ be the differential equation of the $j$-function (see Section \ref{jBackground}). Consider its two-variable version $E_j(x,y)$ given by $$F_j(y, \partial_x y, \partial^2_x y, \partial^3_x y)=0$$ where, as above, $\partial_x = \frac{1}{Dx}\cdot D$. Pila and Tsimerman proved in \cite{Pila-Tsim-Ax-j} that $E_j$ satisfies an Ax-Schanuel inequality of the above form where $\mathcal{P}$ is the collection of all modular polynomials. Hence the above result implies that the set $F_j(y,Dy,D^2y,D^3y)=0$ is strongly minimal and geometrically trivial.  {This was first proven by Freitag and Scanlon in \cite{Freitag-Scanlon}, and for the $j$-function our work recovers their proof modulo some minor differences discussed later in the paper.}

%In fact, the Pila-Tsimerman inequality is the main motivation for this paper.

Thus we get a necessary condition for $E$ to satisfy an Ax-Schanuel inequality of the given form. This is a step towards the solution of the problem described above. In particular, it gives rise to a converse problem: given a one-variable differential equation which is strongly minimal and geometrically trivial, can we say anything about the Ax-Schanuel properties or possibly weaker transcendence properties, such as the Ax-Lindemann-Weierstrass property\footnote{ {As we will see in Section \ref{Part 1}, we only need the Ax-Lindemann-Weierstrass property in the proof of Theorem \ref{thm:intro-AS-SM}.}} considered later, of its two-variable analogue? We briefly discuss this question in Section \ref{remarks}.  {The relation between strong minimality of certain differential equations and transcendence properties of their sets of solutions has been recently studied in \cite{cas-frei-nag,sanz-cas-frei-nag}.}

Further, understanding the structure of strongly minimal sets in a given theory is a central problem in geometric model theory. In $\DCF_0$ there is a nice classification of strongly minimal sets. Namely, they satisfy the Zilber trichotomy, that is, such a set must be either geometrically trivial or non-orthogonal to a Manin kernel\footnote{More precisely, it is non-orthogonal to the Manin kernel $A^\#$ of a simple abelian variety $A$ of $C$-trace zero.} (this is the locally modular non-trivial case) or non-orthogonal to the field of constants which corresponds to the non-locally modular case \cite{Hru-Sok}. Hrushovski \cite{Hru-Juan} also gave a full characterisation of strongly minimal sets of order $1$ proving that such a set is either non-orthogonal to the constants or it is trivial and $\aleph_0$-categorical. However, there is no general classification of trivial strongly minimal sets of higher order and therefore we do not fully understand the nature of those sets. From this point of view the set $J$ defined by the differential equation of $j$ is quite intriguing since it is the first example of a trivial strongly minimal set in $\DCF_0$ which is not $\aleph_0$-categorical. Before Freitag and Scanlon established those properties of $J$ in \cite{Freitag-Scanlon}, it was mainly believed that trivial strongly minimal sets in $\DCF_0$ must be $\aleph_0$-categorical. The reason for this speculation was Hrushovski's aforementioned theorem on order $1$ strongly minimal sets and the lack of counterexamples.

Thus, the classification of strongly minimal sets in $\DCF_0$ can be seen as another source of motivation for the work in this paper, where we show that these two problems (Ax-Schanuel type theorems and geometry of strongly minimal sets) are in fact closely related.

In the second part of the paper we use Ax-Schanuel for the $j$-function to classify strongly minimal sets in $j$-reducts of differentially closed fields (this question was asked by Zilber in a private communication). More precisely, the problem is to classify strongly minimal sets in the reduct $\K:=(K;+, \cdot, E_j)$ of $\mathcal{K}$ where $E_j$ is the two-variable differential equation of the $j$-function described above.  We establish the following dichotomy result as an answer to that question.

\begin{theorem}\label{thm:intro-K-dichotomy}
In $\K$ all strongly minimal sets are geometrically trivial or non-orthogonal to the field of constants $C$ (the latter being definable in $\K$).
\end{theorem}

The proof of this theorem is based not only on Ax-Schanuel, but also on an Existential Closedness statement (see \cite{Aslanyan-Eterovic-Kirby-Diff-EC-j,Aslanyan-adequate-predim}), which asserts that if for a system of equations in $\K$ having a solution does not contradict Ax-Schanuel then it does have a solution.
The Existential Closedness statement is related to the question of adequacy of the Ax-Schanuel inequality for the $j$-function (see \cite[\S 2]{Aslanyan-adequate-predim}). Adequacy means that the Ax-Schanuel inequality governs the geometry of the reduct, hence it is not surprising that it leads to a classification of strongly minimal sets there.

We also study strongly minimal sets in a more basic reduct, namely $\mathcal{K}_C:=(K;+, \cdot, C)$ where $C$ is the field of constants, which is just a pair of algebraically closed fields of characteristic $0$. Actually, this is the first example that we deal with in the second part of the paper. For this reduct we do not have any Ax-Schanuel type statement and we do not need one since it is quite easy to understand definable sets in such a structure. In this case we have the following result.

\begin{theorem}\label{thm:pairs-dichotomy}
All strongly minimal sets in $\mathcal{K}_C$ are non-orthogonal to $C$.
\end{theorem}

 {Although we have not seen this theorem in the literature, we believe it is well known to experts. In fact, as we will see, it easily follows from some basic observations on pairs of algebraically closed fields. Thus, we merely present our proof of Theorem \ref{thm:pairs-dichotomy} as a prelude to the aforementioned classification of strongly minimal sets in $j$-reducts.}

The paper is organised as follows. In Section \ref{jBackground} we give a brief account of the $j$-function. Section \ref{Part 1} is the ``first part'' of the paper where we study strong minimality and geometric triviality of certain definable sets in $\DCF_0$. Section \ref{Part 2}, the ``second part'', is devoted to the classification of strongly minimal sets in $j$-reducts of $\DCF_0$. Appendix \ref{A} contains some preliminaries on strongly minimal sets.

\subsection*{Notation and conventions}

\begin{itemize}

    \item The length of a tuple $\bar{a}$ will be denoted by $| \bar{a}|$.

    \item For a set $A$ and a tuple $\bar{a}\in A^n$ we will often write $\bar{a}\subseteq A$ when the length of $\bar{a}$ is not important.

    \item  {In this paper all fields are of characteristic $0$.}
    
    \item For fields $L\subseteq K$ the transcendence degree of $K$ over $L$ is denoted by $\td(K/L)$ or $\td_LK$. The algebraic locus (Zariski closure) of a tuple $\bar{a} \in K^n$ over $L$ will be denoted by $\Loc_L (\bar{a})$ or $\Loc(\bar{a}/L)$.
    
    \item The algebraic closure of a field $L$ is denoted by $L^{\alg}$.
    
    \item Algebraic varieties defined over an (algebraically closed) field $L$ will be identified with the sets of their $L$-rational points.
    
    \item In a differential field $(K;+, \cdot, D)$ and a non-constant element $x$ the \emph{differentiation with respect to} $x$ is a derivation $\partial_x$ of $K$ defined by $\partial_x : y \mapsto \frac{Dy}{Dx}$.

    \item $\Mr$ stands for Morley rank.

\end{itemize}

\section{Background on the $j$-function}\label{jBackground}

The $j$-function is a modular function for $\SL_2(\mathbb{Z})$, which is defined and analytic on the upper half-plane $\mathbb{H}:=\{ z \in \mathbb{C}: \im(z)>0 \}$. Let $\GL_2^+(\mathbb{Q})$ be the subgroup of $\GL_2(\mathbb{\mathbb{Q}})$ consisting of matrices with positive determinant. This group acts on the upper half-plane via linear fractional transformations. For $g \in \GL_2^+(\mathbb{Q})$ we let $N(g)$ be the determinant of $g$ scaled so that it has relatively prime integral entries. For each positive integer $N$ there is an irreducible polynomial $\Phi_N(X,Y)\in \mathbb{Z}[X,Y]$ such that whenever $g \in \GL_2^+(\mathbb{Q})$ with $N=N(g)$, the function $\Phi_N(j(z),j(gz))$ is identically zero. Conversely, if $\Phi_N(j(z_1),j(z_2))=0$ for some $z_1,z_2 \in \mathbb{H}$ then $z_2=gz_1$ for some $g \in \GL_2^+(\mathbb{Q})$ with $N=N(g)$. The polynomials $\Phi_N$ are called \emph{modular polynomials} (see \cite{Lang-elliptic}). It is well known that $\Phi_1(X,Y)=X-Y$ and all the other modular polynomials are symmetric. %For $w=j(z)$ the image of the $\GL_2^+(\mathbb{Q})$-orbit of $z$ under $j$ is called the \emph{Hecke orbit} of $w$. It obviously consists of the union of solutions of the equations $\Phi_N(w,X)=0,~ N\geq 1$. 

\begin{definition}
Two elements $w_1,w_2 \in \mathbb{C}$ are called \emph{modularly independent} if they do not satisfy any modular relation $\Phi_N(w_1,w_2)=0$. 
\end{definition}
This definition makes sense for arbitrary fields of characteristic zero as the modular polynomials have integer coefficients.

The $j$-function satisfies a third order algebraic differential equation over $\mathbb{Q}$, and none of lower order, i.e. its differential rank over $\mathbb{C}$ is $3$. Namely, $F_j(j,j',j'',j''')=0$ where 
$$F_j(Y_0,Y_1,Y_2,Y_3)=\frac{Y_3}{Y_1}-\frac{3}{2}\left( \frac{Y_2}{Y_1} \right)^2 + \frac{Y_0^2-1968Y_0+2654208}{2Y_0^2(Y_0-1728)^2}\cdot Y_1^2.$$

Thus
$$F_j(Y,Y',Y'',Y''')=SY+R(Y)(Y')^2,$$
where $S$ denotes the \emph{Schwarzian derivative} defined by $$SY = \frac{Y'''}{Y'} - \frac{3}{2} \left( \frac{Y''}{Y'} \right) ^2$$ and $$R(Y)=\frac{Y^2-1968Y+2654208}{2Y^2(Y-1728)^2}$$ is a rational function. All functions $j(gz)$ with $g\in \SL_2(\mathbb{C})$ satisfy this equation and all solutions are of that form (if one wants a solution to be defined on $\mathbb{H}$ then one takes $g\in \SL_2(\mathbb{R})$). See, for example, \cite[Lemma 4.2]{Freitag-Scanlon} or \cite[Lemma 4.1]{Aslanyan-adequate-predim} for a proof.

Here $'$ denotes the derivative of a complex function. When we work in an abstract differential field we will always denote its derivation by $D$ and for an element $a$ in that field $a', a'',\ldots$ will be some other elements and not necessarily the derivatives of $a$.

In an abstract differential field $(K;+,\cdot,D,0,1)$ the differential equation of $j$ is the equation $F_j(y,Dy,D^2y,D^3y)=0$.  {Let $$f_j(x,y):=F_j\left(y,\partial_xy,\partial^2_xy,\partial^3_xy\right).$$ Then $f_j(x,y)=0$ is the two-variable equation of the $j$-function. When $K$ is a field of meromorphic functions (of some variable $t$) then a solution of this equation is a pair of functions $(x(t),y(t))$ where $y(t) = j(gx(t))$ for some $g\in \SL_2(\mathbb{C})$.}

\begin{theorem}[Ax-Schanuel with Derivatives for $j$, {\cite[Theorem 1.3]{Pila-Tsim-Ax-j}}]\label{j-chapter-Ax-for-j}
Let $z_i, j_i \in K \setminus C,~ i=1,\ldots,n,$ be such that $f_j(z_i,j_i)=0.$ If $j_i$'s are pairwise modularly independent then 
\begin{equation*}\label{j-chapter-Ax-ineq}
\td_CC (z_i,\partial_{z_i}j_i,\partial^2_{z_i}j_i : 1\leq i \leq n) \geq 3n+1.
\end{equation*}
\end{theorem}

 {Note that \cite[Theorem 1.3]{Pila-Tsim-Ax-j} is in fact the Ax-Schanuel theorem for several commuting derivations, and the above statement is a special case of that.}

\section{Ax-Schanuel and geometry of strongly minimal sets in $\DCF_0$}\label{Part 1}

\subsection{Setup and main results}\label{j-setup}

Recall that $\mathcal{K}=(K;+,\cdot, D, 0, 1)$ is a differentially closed field with field of constants $C$. We may assume without loss of generality that $\mathcal{K}$ is sufficiently saturated if necessary. Fix an element $t$ with $Dt=1$. Let $E(x,y)$ be the set of solutions of a differential equation $f(x,y)=0$ with constant coefficients.

We give several definitions and then state the main results of the first part of the paper. 

\begin{definition}
Let $\mathcal{P}$ be a non-empty collection of algebraic polynomials $P(X,Y) \in C[X,Y]$. We say two elements $a, b \in K$ are $\mathcal{P}$-\emph{independent} if $P(a,b) \neq 0$ and $P(b,a) \neq 0$ for all $P \in \mathcal{P}$. The $\mathcal{P}$-\emph{orbit} of an element $a \in K$ is the set $\{ b \in K: P(a,b)=0 \mbox{ or } P(b,a)=0 \mbox{ for some } P \in \mathcal{P}\}$ (in analogy to a Hecke orbit). Also, $\mathcal{P}$ is said to be \emph{trivial} if it consists only of the polynomial $X-Y$.
\end{definition}

Recall that $f(x,y)=0$ is the differential equation defining $E$ and let $m:= \ord_Yf(X,Y)$ be the order of $f$ with respect to $Y$.

\begin{definition}
\begin{itemize}[leftmargin=0.45cm]
    \item[]
    
    \item  We say $E(x,y)$ has the  $\mathcal{P}$-ASD property (Ax-Schanuel with Derivatives with respect to $\mathcal{P}$) if the following condition is satisfied.
    \begin{itemize}[leftmargin=0.4cm]
    \item[] 

    Let $x_1,\ldots,x_n,y_1,\ldots,y_n$ be non-constant elements of $K$ with $f(x_i,y_i)=0$ for every $i$. If the $y_i$'s are pairwise $\mathcal{P}$-independent then 
    \begin{equation}\label{general-Ax-ineq}
    \td_CC\left(x_1, y_1, \partial_{x_1}y_1, \ldots, \partial^{m-1}_{x_1}y_1, \ldots, x_n, y_n, \partial_{x_n}y_n, \ldots, \partial^{m-1}_{x_n}y_n\right) \geq mn+1.
    \end{equation}
    \end{itemize}
    
    \item  We say $E$ has the  $\mathcal{P}$-ALWD property (Ax-Lindemann-Weierstrass with Derivatives with respect to $\mathcal{P}$)  if the inequality \eqref{general-Ax-ineq} is satisfied under the additional assumption $\td_CC(\bar{x})=1$.
 \end{itemize}
\end{definition}

\begin{remark}
 {In other words, the $\mathcal{P}$-ALWD property states that if $x,y_1,\ldots,y_n\in K$ are non-constant elements such that $x_1,\ldots,x_n \in C(x)^{\alg}\setminus C$ and $f(x_i,y_i)=0$ and $y_1,\ldots,y_n$ are pairwise $\mathcal{P}$-independent, then the $mn$ elements $$y_1, \partial_{x_1}y_1, \ldots, \partial^{m-1}_{x_1}y_1, \ldots, y_n, \partial_{x_n}y_n, \ldots, \partial^{m-1}_{x_n}y_n$$ are algebraically independent over $C(x)$.}
\end{remark}

The  $\mathcal{P}$-ASD property can be reformulated as follows: for any non-constant solutions $(x_i,y_i)$ of $E$ the transcendence degree in \eqref{general-Ax-ineq} is strictly larger than $m$ times the number of distinct $\mathcal{P}$-orbits of $y_i$'s. Note that \eqref{general-Ax-ineq} is motivated by the known examples of Ax-Schanuel inequalities  {restricted to ordinary differential fields (as opposed to fields with several derivations) \cite{Ax,Pila-Tsim-Ax-j,Aslanyan-linear}. Similarly, $\mathcal{P}$-ALWD is the abstract version of Pila's Modular Ax-Lindemann-Weierstrass with Derivatives theorem, again for ordinary differential fields. As we explain in Subsection \ref{j}, the differential equation of the $j$-function has both the $\mathcal{P}$-ASD and $\mathcal{P}$-ALWD properties, and the latter is equivalent to a special case of \cite[Theorem 1.1]{Pila-modular}.}

\begin{remark}\label{PvsPredim}
Having the  $\mathcal{P}$-ASD property for a given equation $E$ may force $\mathcal{P}$ to be ``closed'' in some sense. Firstly, $X-Y$ (or a polynomial divisible by $X-Y$) must be in $\mathcal{P}$. Secondly, if $P_1, P_2 \in \mathcal{P}$ then $P_1(y_1,y_2)=0,~ P_2(y_2,y_3)=0$ impose a relation on $y_1$ and $y_3$ given by $Q(y_1,y_3)=0$ for some polynomial $Q$. Then the  $\mathcal{P}$-ASD property may fail if $Q\notin \mathcal{P}$. %one allows a relation $Q(y_i,y_j)=0$ between $y_i$ and $y_j$ (although one requires $P_l(y_i,y_j)\neq 0,~ l=1,2$). 
In that case one has to add $Q$ to $\mathcal{P}$ in order to allow the possibility of an Ax-Shcanuel property with respect to $\mathcal{P}$. 

Similar conditions on $\mathcal{P}$ are required in order for $\mathcal{P}$-independence to define a dimension function of a pregeometry (number of distinct $\mathcal{P}$-orbits), which would imply that the  $\mathcal{P}$-ASD property is a predimension inequality. Note that the collection of modular polynomials has all those properties. However, the shape of $\mathcal{P}$ is not important for our results since we assume that a given equation $E$ has the  $\mathcal{P}$-ASD property.
\end{remark}

\begin{definition}
A $\mathcal{P}$-\emph{special variety} in $K^n$ for some $n$ is an irreducible  component (over  $C$) of a Zariski closed set in $K^n$ defined by a finite collection of equations of the form $P_{ik}(y_i,y_k)=0$ for some $P_{ik} \in \mathcal{P}$. For a subfield $L \subseteq K$ a \emph{weakly} $\mathcal{P}$-\emph{special variety over} $L$ is an irreducible  component over  $L^{\alg}$ of a Zariski closed set in $K^n$ defined by a finite collection of equations of the form $P_{ik}(y_i,y_k)=0$ and $y_i=a$ for some $P_{ik} \in \mathcal{P}$ and $a \in L^{\alg}$.
For a definable set $V$, a (weakly) $\mathcal{P}$-\emph{special subvariety} (over $L$) of $V$ is an intersection of $V$ with a (weakly) $\mathcal{P}$-special variety (over $L$).

A $\mathcal{P}$-special variety $S$ may have a constant coordinate defined by an equation $P(y_i,y_i)=0$ for some $P\in \mathcal{P}$. If no coordinate is constant on $S$ then it is said to be \emph{strongly} $\mathcal{P}$-special.
\end{definition}

Let $C_0 \subseteq C$ be the subfield of $C$ generated by the coefficients of $f$ and let $K_0:=C_0 \langle t \rangle = C_0(t)$ be the differential subfield of $K$ generated by $C_0$ and $t$, where $t$ is an element with $Dt=1$. We fix $K_0$ and work over it; in other words we expand our language with new constant symbols for elements (generators) of $K_0$.

Now we can formulate one of our main results (see Appendix \ref{A} for definitions of geometric triviality and strict disintegratedness).

\begin{theorem}\label{main-thm}
Assume $E(x,y)$ satisfies the  $\mathcal{P}$-ALWD property for some $\mathcal{P}$. Assume further that the differential polynomial $g(Y):=f(t,Y)$ is absolutely irreducible.
Then 
\begin{itemize}
\item $U:=\{y: g(y)=0 \wedge Dy \neq 0 \}$ is strongly minimal with trivial geometry.
\item If, in addition, $\mathcal{P}$ is trivial then $U$ is strictly disintegrated and hence it has $\aleph_0$-categorical induced structure.
\item All definable subsets of $U^n$ over a differential field $L \supseteq K_0$ are Boolean combinations of weakly $\mathcal{P}$-special subvarieties over $L$. 
\end{itemize}
\end{theorem}

\begin{remark}
If the polynomials from $\mathcal{P}$ have rational coefficients then $\mathcal{P}$-special varieties are defined over $\mathbb{Q}^{\alg}$. Furthermore, if $E$ satisfies the  $\mathcal{P}$-ASD property then $U \cap C(t)^{\alg} = \emptyset$ and so $\mathcal{P}$-special subvarieties of $U$ over $C(t)$ are merely $\mathcal{P}$-special subvarieties (over $C$).
\end{remark}
%We first prove this theorem for $s=t$ and then deduce the result for an arbitrary $s$ from this. 
As the reader may guess and as we will see in the proof, this theorem holds under weaker assumptions on $E$. Namely, it is enough to require that \eqref{general-Ax-ineq} hold for $x_1=\ldots=x_n=t$ which is a weak form of the  $\mathcal{P}$-ALWD property. %However, we prefer the given formulation of Theorem \ref{main-thm} since the main object of our interest is the Ax-Schanuel inequality (for $E$).

Further, we deduce from Theorem \ref{main-thm} that if $E$ has some special form, then all fibres $E(s,y)$ for a non-constant $s \in K$ have the above properties (over $C_0\langle s \rangle)$.

\begin{theorem}\label{cor}
Let $E(x,y)$ be defined by $R(x, y, \partial_x y, \ldots, \partial_x^m y)=0$ where $R(X,\bar{Y})$ is an algebraic polynomial over $C$, irreducible over $C(X)^{\alg}$ (as a polynomial of $\bar{Y}$). Assume $E(x,y)$ satisfies the  $\mathcal{P}$-ALWD property for some $\mathcal{P}$ and let $s \in K$ be a non-constant element. Then 
\begin{itemize}
\item $U_s:=\{y: E(s,y) \wedge Dy \neq 0 \}$ is strongly minimal with trivial geometry.
\item If, in addition, $\mathcal{P}$ is trivial then any distinct non-algebraic (over $C_0 \langle s \rangle)$ elements are independent and $U_s$ is $\aleph_0$-categorical.
\item All definable subsets of $U_s^n$ over a differential field $L \supseteq C_0\langle s \rangle$ are Boolean combinations of weakly $\mathcal{P}$-special subvarieties over $L$.
\end{itemize}
\end{theorem}

\begin{remark}
Since $U_s \cap C = \emptyset$, in Theorems \ref{main-thm} and \ref{cor} the induced structure on $U_s^n$ is actually given by \emph{strongly $\mathcal{P}$-special} subvarieties (over $L$), which means that we do not allow any equation of the form $y_i=c$ for $c$ a constant. In particular we also need to exclude equations of the form $P(y_i,y_i)=0$ for $P \in \mathcal{P}$.
\end{remark}

We also prove a generalisation of Theorem \ref{main-thm}.

\begin{theorem}\label{main-gen}
Assume $E(x,y)$ satisfies the  $\mathcal{P}$-ALWD property and let $p(Y)\in C(t)[Y] \setminus C,~ q(Y)\in C[Y] \setminus C$ be such that the differential polynomial $f(p(Y),q(Y))$ is absolutely irreducible. Then the set $$U_{p,q}:=\{y: E(p(y),q(y)) \wedge y\notin C \}$$ is strongly minimal and geometrically trivial.
\end{theorem}

%definable in the language of rings over $C$, that is, for each $n, l > 0$ the set $\{(x_1,\ldots,x_l)\in K^l: d(\bar{x})\geq n\}$ is type-definable in the language of rings over $C$.

As an application of Theorem \ref{main-thm} we obtain a result on the differential equation of the $j$-function which was first established by Freitag and Scanlon in \cite{Freitag-Scanlon}. %To be more precise, let $F_j(j,j',j'',j''')=0$ be the algebraic differential equation satisfied by the modular $j$-function.

\begin{theorem}[{\cite[Theorems 4.5 and 4.7]{Freitag-Scanlon}}]\label{j-trivial}
The set $J \subseteq K$ defined by the equation $F_j(y,Dy,D^2y,D^3y)=0$ is strongly minimal with trivial geometry. Furthermore, $J$ is not $\aleph_0$-categorical.
\end{theorem}

Strong minimality and geometric triviality of $J$ follow directly from Theorem \ref{main-thm} combined with the Ax-Schanuel theorem for $j$ or the Ax-Lindemann-Weierstrass with Derivatives theorem (see \cite[Theorem 1.1]{Pila-modular} and Section \ref{j} of this paper). Of course, the ``furthermore'' clause does not follow from Theorem \ref{main-thm} but it is not difficult to prove. Theorem \ref{main-thm} also gives a characterisation of the induced structure on the Cartesian powers of $J$. Again, that result can be found in \cite{Freitag-Scanlon}.% in a more general form.

 {The proof of Theorem \ref{j-trivial} by Freitag and Scanlon is also based on Pila's modular Ax-Lindemann-Weierstrass with Derivatives theorem, and our proof in that case is quite similar to theirs. There are some minor differences though which are discussed in Section \ref{j}.}

\subsection{Proofs of the main results}\label{j-proofs}

\subsubsection*{Proof of Theorem \ref{main-thm}}

Taking $x_1=\ldots=x_n=t$ in the  $\mathcal{P}$-ALWD property we get the following  property for $U$ which in fact is enough to prove Theorem \ref{main-thm}.

\begin{lemma}\label{Ax-Lindemann}
The $\mathcal{P}$-ALWD property implies that for any pairwise $\mathcal{P}$-independent elements $u_1, \ldots, u_n \in U$ the elements $\bar{u}, D\bar{u}, \ldots, D^{m-1}\bar{u}$ are algebraically independent over $C(t)$, and hence over $K_0$, where $\bar{u}:=(u_1,\ldots,u_n)$.
\end{lemma}

%\begin{proof}[Proof of Theorem \ref{main-thm}]
We show that every definable (possibly with parameters) subset of $U$ is either finite or co-finite.  {Since $U$ is defined over $K_0$, by stable embedding for any definable subset $V\subseteq U$ there is a finite subset $A=\{a_1, \ldots, a_n\} \subseteq U$ such that $V$ is defined over $K_0\cup A$.} It suffices to show that there is a unique non-algebraic type over $K_0\cup A$ realised in $U$, i.e. for any $u_1, u_2 \in U \setminus \acl(K_0 \cup A)$ we have $$\tp(u_1/K_0 \cup A)=\tp(u_2/K_0\cup A).$$ Let $ u \in U \setminus \acl(K_0\cup A)$. We know that $$\acl(K_0\cup A)=K_0\langle A \rangle^{\alg} = K_0(\bar{a}, D\bar{a}, \ldots, D^{m-1}\bar{a})^{\alg}.$$ Since $u \notin K_0\langle A \rangle^{\alg}$, $u$ is transcendental over $K_0(A)$ and hence it is $\mathcal{P}$-independent from each $a_i$. We may assume without loss of generality that $a_i$'s are pairwise $\mathcal{P}$-independent (otherwise we could replace $A$ by a maximal pairwise $\mathcal{P}$-independent subset). Applying Lemma \ref{Ax-Lindemann} to $\bar{a},u$, we deduce that $u,Du,\ldots, D^{m-1}u$ are algebraically independent over $K_0\langle A\rangle$. Hence $\tp(u/K_0\cup A)$ is determined uniquely (axiomatised) by the set of formulae
$$\{ g(y)=0 \} \cup \{ h(y) \neq 0: h(Y) \in K_0\langle A \rangle\{ Y\},~ \ord (h) < m \}$$
In other words, $g(Y)$ is the minimal differential polynomial of $u$ over $K_0\langle A \rangle$ (recall that $g$ is absolutely irreducible and hence it is irreducible over any field).

Thus, $U$ is strongly minimal. {A similar} argument shows also that if $ A \subseteq U$ is a (finite) subset and $ u \in U \cap \acl(K_0A)$ then there is $a \in A$ such that $u \in \acl(K_0a)$. This proves that $U$ is geometrically trivial.  {Note that when $U$ is defined over the constants and has order $>1$, which is the case for the $j$-function, strong minimality automatically implies geometric triviality (see, for example, \cite[Proposition 5.8]{cas-frei-nag}). However, that fact is based on the trichotomy theorem in differentially closed fields, while the above argument shows that triviality is a direct consequence of the $\mathcal{P}$-ALWD property.}

If $\mathcal{P}$ is trivial then distinct elements of $U$ are independent, hence $U$ is strictly disintegrated.

The last part of Theorem \ref{main-thm} follows from the following lemma.
%\end{proof}

\begin{lemma}\label{Induced-structure}
Every irreducible relatively Kolchin closed subset of $U^n$ over $C(t)$ is a $\mathcal{P}$-special subvariety of $U^n$.
\end{lemma}

\begin{proof}
Let $V \subseteq U^n$ be an irreducible relatively closed subset, i.e. it is the intersection of $U^n$ with an irreducible Kolchin closed set in $K^n$. Pick a generic point $\bar{v}=(v_1, \ldots, v_n) \in V$ and let $W \subseteq K^n$ be the Zariski closure of $\bar{v}$ over $C$. Let $d:=\dim W$ and assume $v_1, \ldots, v_d$ are algebraically independent over $C$. Then $v_i \in C(v_1, \ldots, v_d)^{\alg}$ for each $i=d+1, \ldots, n$. By Lemma \ref{Ax-Lindemann} each $v_i$ with $i>d$ must be in a $\mathcal{P}$-relation with some $v_{k_i}$ with $k_i \leq d$. Let $P_i(v_i, v_{k_i})=0$ for $i>d$. The algebraic variety defined by the equations $P_i(y_i, y_{k_i})=0,~ i=d+1, \ldots, n,$ has dimension $d$ and contains $W$. Therefore $W$ is a component of that variety and so it is a $\mathcal{P}$-special variety.

We claim that $W \cap U^n = V$. Since $v_1, \ldots, v_d \in U$ are algebraically independent over $C$, by Lemma \ref{Ax-Lindemann} $\bar{v}, D\bar{v}, \ldots, D^{m-1}\bar{v}$ are algebraically independent over $C(t)$. Moreover, the (differential) type of each $v_i,~ i>d,$ over $v_1, \ldots, v_d$ is determined uniquely by an irreducible algebraic equation. Therefore $\tp(\bar{v}/C(t))$ is axiomatised by formulas stating that $\bar{v}$ is Zariski generic in $W$ and belongs to $U^n$. In other words, $\bar{v}$ is Kolchin generic in $W \cap U^n$. Now $V$ and $W\cap U^n$ are both equal to the Kolchin closure of $\bar{v}$ inside $U^n$ and hence they are equal.
\end{proof}

Thus, definable subsets of $U^n$ over $C(t)$ are Boolean combinations of special subvarieties. Now let $L \subseteq K$ be an arbitrary differential subfield over which $U$ is defined. Then definable subsets of $U^n$ over $L$ can be defined with parameters from $\tilde{L} = K_0 \cup (U \cap L^{\alg})$ (see Appendix \ref{A}). Then Lemma \ref{Induced-structure} implies that irreducible Kolchin closed subsets of $U^n$ defined over $\tilde{L}$ are weakly $\mathcal{P}$-special subvarieties of $U^n$ over $L$.

%components of algebraic varieties defined by equations of the form $P(x_i,x_k)=0$ and $x_i=l$ for some $l \in \tilde{L}$.

Finally, note that since $U$ does not contain any algebraic elements over $C(t)$, the type of any element $u \in U$ over $C(t)$ is isolated by the formula $f(t,y)=0 \wedge Dy \neq 0$.

\subsubsection*{Proof of Theorem \ref{main-gen}}

We argue as above and show that for a finite set $$A=\{a_1,\ldots,a_n\} \subseteq U_{p,q}$$ there is a unique non-algebraic type over $K_0\langle A \rangle$ realised in $U_{p,q}$. Here we will use the full $\mathcal{P}$-ALWD property. 

%It is not difficult to observe that assuming there are two $d$-independent elements in $K$ and the set $\{ (y_1,y_2)\in K^2: d(y_1,y_2)=2\}$ is type-definable over $C$, any two algebraically independent elements over $C$ are $d$-independent. Therefore 
If $u \in U_{p,q} \setminus (K_0\langle A \rangle)^{\alg}$ then $q(u)$ is transcendental over $K_0(A)$ and so $q(u)$ is $\mathcal{P}$-independent from each $q(a_i)$.  Moreover, we may assume $\{ q(a_1),\ldots,q(a_n) \} $ is $\mathcal{P}$-independent. Then by the  $\mathcal{P}$-ALWD property 
$$\td_CC\left(p(u), q(u), \ldots, \partial^{m-1}_{p(u)}q(u), p(a_i), q(a_i), \ldots, \partial^{m-1}_{p(a_i)}q(a_i)\right)_{i=1,\ldots,n} \geq m(n+1)+1.$$
But then $$\td_CC\left(t,u,Du,\ldots,D^{m-1}u,a_i,Da_i,\ldots,D^{m-1}a_i\right)_{i=1,\ldots,n} \geq m(n+1)+1,$$
and hence $u,Du,\ldots,D^{m-1}u$ are algebraically independent over $K_0\langle A \rangle$. This determines the type $\tp(u/K_0A)$ uniquely as required. It also shows triviality of the geometry.

\subsubsection*{Proof of Theorem \ref{cor}}

Consider the differentially closed field $\mathcal{K}_s := (K;+,\cdot, \partial_s,0,1)$. 
The given form of the differential equation $E$ implies that $U_s$ is defined over $C_0(s)$ in $\mathcal{K}_s$. However, in general it may not be defined over $C_0(s)$ in $\mathcal{K}$, it is defined over $C_0 \langle s \rangle = C_0(s,Ds,D^2s,\ldots)$. As $s\notin C$, it is transcendental over $C$ and so $R(s,\bar{Y})$ is irreducible over $C(s)^{\alg}$. Therefore $R(s, Y, \partial_s Y, \ldots, \partial_s^m Y)$ is absolutely irreducible. Since $\partial_ss=1$, we know by Theorem \ref{main-thm} that $U_s$ is strongly minimal in $\mathcal{K}_s$. On the other hand, the derivations $\partial_s$ and $D$ are inter-definable (with parameters) and so a set is definable in $\mathcal{K}$ if and only if it is definable in $\mathcal{K}_s$, possibly with different parameters. This implies that every definable subset of $U_s$ in $\mathcal{K}$ is either finite or co-finite, hence it is strongly minimal.

Further, Theorem \ref{main-thm} implies that $U_s$ is geometrically trivial over $C_0(s)$ in $\mathcal{K}_s$. By Theorem \ref{geometric-triviality}, $U_s$ is also geometrically trivial over $C_0 \langle s \rangle$ in $\mathcal{K}_s$. On the other hand, for any subset $A \subseteq U_s$ the algebraic closure of $C_0 \langle s \rangle \cup A$ in $\mathcal{K}$ is the same as in $\mathcal{K}_s$. This implies geometric triviality of $U_s$ in $\mathcal{K}_s$.

The same argument (along with the remark after Theorem \ref{geometric-triviality}) shows that the second and the third parts of Corollary \ref{cor} hold as well.

\subsection{Application to the $j$-function}\label{j}

Recall that the differential equation of the $j$-function is of the form
\begin{equation}\label{dif-eq-j}
Sy+R(y)(Dy)^2 = 0 
\end{equation}
where $S$ denotes the Schwarzian derivative and $R$ is a rational function. Let $J$ be the set defined by \eqref{dif-eq-j}. Note that $F_j(y,Dy,D^2y,D^3y)=Sy+R(y)(Dy)^2$ is a differential rational function, and not a polynomial. In particular, constant elements do not satisfy \eqref{dif-eq-j} for $Sy$ is not defined for a constant $y$. We can multiply our equation through by a common denominator and make it into a polynomial equation $F_j^*(y,Dy,D^2y,D^3y)=0$ with
\begin{equation}\label{j-eq-poly}
F_j^*(y,Dy,D^2y,D^3y):= q(y)DyD^3y-\frac{3}{2}q(y)(D^2y)^2+p(y)(Dy)^4,
\end{equation}
where $p$ and $q$ are respectively the numerator and the denominator of $R$. Let $J^*$ be the set defined by \eqref{j-eq-poly}. It is not strongly minimal since $C$ is a definable subset. However, as we will see shortly, $J=J^*\setminus C$ is strongly minimal and $\Mr(J^*)=1,~ \MD(J^*)=2$. Thus, whenever we speak of the formula $F_j(y,Dy,D^2y,D^3y)=0$ (which, strictly speaking, is not a formula in the language of differential rings), we mean the formula $F_j^*(y,Dy,D^2y,D^3y)=0 \wedge Dy \neq 0$. %Note also that $J^*$ is Kolchin closed while $J$ is Kolchin constructible. 

Let $\Phi:=\{ \Phi_N(X,Y): N>0\}$ be the collection of modular polynomials. Then two elements are modularly independent if and only if they are $\Phi$-independent. For an element $a \in K$ its \emph{Hecke orbit} is its $\Phi$-orbit. 

Consider the two-variable analogue of the equation \eqref{j-eq-poly}:
\begin{equation}\label{Two-var-eq}
F^*_j(y,\partial_xy,\partial^2_xy,\partial^3_xy)=0.
\end{equation}

Theorem \ref{j-chapter-Ax-for-j} states that \eqref{Two-var-eq} has the $\Phi$-ASD property and hence also the $\Phi$-ALWD property. Thus, as a consequence of Theorems \ref{main-thm} and \ref{j-chapter-Ax-for-j} we get strong minimality and geometric triviality of $J$ (note that $F_j^*(Y_0,Y_1,Y_2,Y_3)$ is absolutely irreducible for it depends linearly on $Y_3$). This was first established by Freitag and Scanlon in \cite{Freitag-Scanlon}. Moreover, Theorem \ref{cor} shows that all non-constant fibres of \eqref{Two-var-eq} are strongly minimal and geometrically trivial (after removing constant points) and the induced structure on the Cartesian powers of those fibres is given by (strongly) special subvarieties. It is proven in \cite{Freitag-Scanlon} that the sets $F_j(y,Dy,D^2y,D^3y)=a$ have the same properties for any $a$.

 {Instead of using the Ax-Schanuel theorem for $j$, we can also deduce the $\Phi$-ALWD property from Pila's modular Ax-Lindemann-Weierstrass with Derivatives theorem  by employing Seidenberg's embedding theorem. Indeed, \cite[Theorem 1.1]{Pila-modular} states that if $\mathbb{C}(W)$ is an algebraic function field and $a_1,\ldots,a_n \in \mathbb{C}(W)$ take values in $\mathbb{H}$ at some point of $W$ and have distinct $\GL_2^+(\mathbb{Q})$-orbits then the functions $$j(a_1),\ldots,j(a_n),j'(a_1),\ldots,j'(a_n),j''(a_1),\ldots,j''(a_n)$$ are algebraically independent over $\mathbb{C}(W)$. When $W$ is a curve, i.e. $\td(\mathbb{C}(W)/C)=1$, this theorem is equivalent to the $\Phi$-ALWD property for the equation $f_j(x,y)=0$. The general form of \cite[Theorem 1.1]{Pila-modular} for arbitrary $W$ corresponds to an Ax-Lindemann-Weierstrass with Derivatives statement in fields with several derivations, which is a special case of the Ax-Schanuel with Derivatives theorem for fields with several derivations.}

 {As it was already mentioned, Freitag and Scanlon also use Pila's Ax-Lindemann-Weierstrass with Derivatives theorem to prove that $J$ is strongly minimal. Therefore, our proof in the case of $j$ is equivalent to their proof, a difference being the model theoretic tools used in the application of Pila's theorem. For example, Freitag and Scanlon use Shelah's reflection principle while we use stable embedding.}

\begin{remark}
To complete the proof of Theorem \ref{j-trivial}, that is, to show that $J$ is not $\aleph_0$-categorical, one argues as follows (see \cite[Theorem 4.7]{Freitag-Scanlon}). The Hecke orbit of an element $ j \in J$ is infinite and is contained in $J$. Therefore $J$ realises infinitely many algebraic types over an arbitrary element $j\in J$ and hence it is not $\aleph_0$-categorical. 
\end{remark}

\subsection{Some remarks}\label{remarks}

An interesting question is whether there are differential equations with the  $\mathcal{P}$-ASD property with trivial $\mathcal{P}$. As we showed here, if $E(x,y)$ has such a property then the corresponding $U$ (and other fibres too) must be strongly minimal and strictly disintegrated. There are quite a few examples of this kind of strongly minimal sets in $\DCF_0$. The two-variable versions of those equations will be natural candidates of equations with the required  $\mathcal{P}$-ASD property. Note that since in our proofs we only used the  $\mathcal{P}$-ALWD  property, it would be more reasonable to expect those equations to satisfy the  $\mathcal{P}$-ALWD property for trivial $\mathcal{P}$. However, as it was mentioned earlier, we are mainly interested in the  $\mathcal{P}$-ASD property of differential equations.

For example, the geometry of the sets of the form $Dy=r(y)$, where $r$ is a rational function over $C$, is well understood. The nature of the geometry is determined by the partial fraction decomposition of $1/r$. As an example consider the equation
\begin{equation}\label{strict-dis-eq}
Dy=\frac{y}{1+y}.
\end{equation}
It defines a strictly disintegrated strongly minimal set (see \cite[Corollary 6.3]{Mar-dif}). The two variable analogue of this equation is
\begin{equation}\label{strict-dis-eq-2}
\partial_x y=\frac{y}{1+y}.
\end{equation}
But this is equivalent to the equation $\frac{Dy}{y} = D(x-y)$. Denoting $x-y$ by $z$ we get the exponential differential equation $Dy=yDz$. It is easy to deduce from this that \eqref{strict-dis-eq-2} does not satisfy the  $\mathcal{P}$-ASD property for any $\mathcal{P}$ (it satisfies a version of the original exponential Ax-Schanuel inequality though). Indeed, the fibre of \eqref{strict-dis-eq-2} above $x=t$ is of trivial type but the section by $x=t+y$ is non-orthogonal to $C$. So according to Theorem \ref{main-gen} the equation \eqref{strict-dis-eq-2} does not satisfy any  $\mathcal{P}$-ALWD property, let alone a $\mathcal{P}$-ASD property. %The fibre of $y'=yx'$ above $x=t$ is non-orthogonal to $C$ but the section by $x=t-y$ is of trivial type. 
Clearly, all the sets $Dy=r(y)$ can be treated in the same manner and hence they are not appropriate for our purpose. Thus, one needs to look at the behaviour of all the sets $E(p(y),q(y))$, and if they happen to be trivial strongly minimal sets then one can hope for a  $\mathcal{P}$-ASD inequality, or at least a $\mathcal{P}$-ALWD property.

The classical Painlev\'{e} equations with generic parameters define strongly minimal and strictly disintegrated sets as well. For example, let us consider the first Painlev\'{e} equation $D^2y=6y^2+t$.  {Strong minimality of this equation was established by Kolchin (see \cite[Theorem 5.18]{Mar-dif}), and algebraic independence of solutions was proven by Nishioka in \cite{Nishioka}.} We consider its two-variable version
\begin{equation}\label{Painleve1}
\partial_x^2y=6y^2+x.
\end{equation}
The goal is to find an Ax-Schanuel inequality for this equation. Observe that \eqref{Painleve1} does not satisfy the  $\mathcal{P}$-ASD property (nor $\mathcal{P}$-ALWD) with trivial $\mathcal{P}$. Indeed, if $\zeta$ is a fifth root of unity then the transformation $x \mapsto \zeta^2 x,~ y \mapsto \zeta y$ sends a solution of \eqref{Painleve1} to another solution. If one believes these are the only relations between solutions of the above equation, then one can pose the following conjecture.

\begin{conjecture}[Ax-Schanuel for the first Painlev\'{e} equation]
If $(x_i,y_i),~ i=1, \ldots, n,$ are solutions to the equation \eqref{Painleve1} and $(x_i/ x_j)^5\neq 1$ for $i\neq j$ then
$$\td(\bar{x},\bar{y},\partial_{\bar{x}}\bar{y}) \geq 2n+1.$$
\end{conjecture}

One could in fact replace $x$'s with $y$'s in the condition $(x_i/ x_j)^5\neq 1$ as those are equivalent. Hence the above conjecture states that \eqref{Painleve1} has the $\{ X^5-Y^5 \}$-ASD property. %Proving the $\{ X^5-Y^5 \}$-ALWD property would give an interesting special case of this conjecture.

Nagloo and Pillay showed in \cite{Nagloo-Pillay} that the other generic Painlev\'{e} equations define strictly disintegrated strongly minimal sets as well. So we can analyse relations between solutions of their two-variable analogues and ask similar questions for them too.

 {In a recent work Casale, Freitag and Nagloo established some functional transcendence results for the uniformising functions of  genus zero Fuchsian groups. In particular, they showed that for a Fuchsian group $\Gamma$ of genus $0$ the differential equation of the uniformiser $j_{\Gamma}$ (which is of the form \eqref{dif-eq-j} where $R$ is a rational function depending on $j_{\Gamma}$) defines a strongly minimal and geometrically trivial set, and it is $\aleph_0$-categorical if and only if the group is non-arithmetic (see \cite[Theorems 2.12, 2.13 and 2.14]{cas-frei-nag}). They also proved the Ax-Lindemann-Weierstrass with Derivatives property for $j_{\Gamma}$ which, in our terminology, is a $\mathcal{P}$-ALWD property where $\mathcal{P}$ consists of the $\Gamma$-special polynomials (see \cite[Theorem 2.16]{cas-frei-nag}). For many non-arithmetic groups $\Gamma$ there are no $\Gamma$-special polynomials, hence $j_{\Gamma}$ satisfies the $\mathcal{P}$-ALWD property with trivial $\mathcal{P}$. Therefore these functions are good candidates for the $\mathcal{P}$-ASD property with trivial $\mathcal{P}$.}

 {Note that while Pila's proof of the ALWD theorem for the $j$-function is based on the theory of o-minimality, the approach of \cite{cas-frei-nag} is different. It does not employ o-minimality and instead makes heavy use of differential algebra and model theory of differential fields. In particular, strong minimality and geometric triviality of the differential equation of $j_{\Gamma}$ is proven first which is then used in the proof of ALWD.}

 {Let us also note that some of these results were refined and generalised by Bl\'asquez-Sanz, Casale, Freitag and Nagloo in \cite{sanz-cas-frei-nag}.}

\section{Strongly minimal sets in $j$-reducts of $\DCF_0$}\label{Part 2}

First we study strongly minimal sets in pairs of algebraically closed fields. It will serve as a simple example of the methods that we are going to use in $j$-reducts.

\subsection{Pairs of algebraically closed fields}

Model theory of pairs of algebraically closed fields has been well studied (see, for example, \cite{vandenDries-pairs-of-fields,Keisler-pairs-of-fields}). Therefore, most of the results of this subsection are known.  {Even though we were not able to find some statements in the existing literature (we give references for those that we could find), we still believe they are known to experts. In any case, the proofs presented here have been obtained independently as a special case of the more delicate analysis of types and strongly minimal sets in $j$-reducts studied in the following subsections.}

Let $\mathcal{K}_C:=(K ;+, \cdot, C )$ be an algebraically closed field of characteristic $0$ with a distinguished algebraically closed subfield $C$ ($C$ is a unary predicate in the language). It is easy to prove that this structure is $\omega$-stable of Morley rank $\omega$ (cf. \cite[Proposition 4.1]{Aslanyan-def-deriv}). We assume $\mathcal{K}_C$ is sufficiently saturated. 

Let $\bar{a}\in K^m$ and $b\in K$. 

\begin{lemma}[cf. {\cite[Proposition 4.2]{vandenDries-pairs-of-fields}}]\label{lem:pairs-MR-fin-alg}
$\Mr(b/\bar{a})<\omega$ iff $b \in C(\bar{a})^{\alg}$.
\end{lemma}
\begin{proof}
If $b$ is transcendental over $C(\bar{a})$ then for any $b' \notin C(\bar{a})^{\alg}$ there is a field automorphism of $K$ fixing $C(\bar{a})$ pointwise and mapping $b$ to $b'$. In particular, it is an automorphism of $\mathcal{K}_C$ and so $\tp(b/\bar{a}) = \tp(b'/\bar{a})$, and this type is the generic type over $\bar{a}$.
\end{proof}

Now let $b\in C(\bar{a})^{\alg}$. Then for some polynomial $p$ the equality $p(\bar{a}, \bar{c}, b) = 0$ holds for some finite tuple $\bar{c}\in C^l$. Let $W:=\Loc_{\mathbb{Q}(\bar{a})} (\bar{c})\subseteq K^l$ be the algebraic locus (Zariski closure) of $\bar{c}$ over $\mathbb{Q}(\bar{a})$. For every proper subvariety $U \subsetneq W$ defined over $\bar{a}$ consider the formula
\begin{equation}\label{eq-exist}
   \varphi_U(y) = \exists \bar{x} (\bar{x} \in C^l \cap (W\setminus U) \wedge p(\bar{a},\bar{x}, y) = 0). 
\end{equation}
Note that $\varphi_U$ implicitly depends on $\bar{a}$ and $b$, and it will be clear from the context what $\bar{a}$ and $b$ are.

Notice that for every $U\subsetneq W$ the formula $\varphi_U(b)$ holds. Observe also that the set $C^l \cap (W(K)\setminus U(K))$, being a subset of $C^l$, is actually definable with parameters from $C$. This follows from the stable embeddedness of $C$ in $K$.

\begin{proposition}[cf. {\cite[Lemma 3.2]{vandenDries-pairs-of-fields}}]\label{prop:pairs-type}
If $b\in C(\bar{a})^{\alg}$, then the collection of all formulas $\varphi_U(y)$ determines (axiomatises) $\tp(b/\bar{a})$.
\end{proposition}
\begin{proof}
Assume $b' \models \varphi_U(y)$ for all $U\subsetneq W$. The collection of formulas $$ \{ \bar{x} \in C^l \cap (W\setminus U) \wedge p(\bar{a},\bar{x}, b') = 0: U\subsetneq W \} $$ (over $\bar{a}, b'$) is finitely satisfiable so it has a realisation $\bar{c}'$. Evidently $\bar{c}'$ is generic in $W$ over $\bar{a}$. Therefore there is an automorphism $\pi$ of $C(\bar{a})$ which fixes $\bar{a}$ pointwise, fixes $C$ setwise and sends $\bar{c}$ to $\bar{c}'$. This automorphism can be extended to an automorphism of $\mathcal{K}_C$ which sends $b$ to $b'$.
\end{proof}

\begin{remark}
This shows, in particular, that the first-order theory\footnote{This theory is axiomatised by axiom schemes stating that $K$ is an algebraically closed field of characteristic $0$ and $C$ is an algebraically closed subfield.} of $\mathcal{K}_C$ is \emph{near model complete}, that is, every formula is equivalent to a Boolean combination of existential formulas modulo that theory. One can also show that it is not model complete. Indeed, pick three algebraically independent elements $a, b, x$ over $\mathbb{Q}$ and set $y:=ax+b$. Let $$C_0:=\mathbb{Q}^{\alg},~ C_1:=\mathbb{Q}(a, b)^{\alg},~ K_0:=\mathbb{Q}(x,y)^{\alg},~ K_1:=\mathbb{Q}(a,b,x)^{\alg}.$$
Then $K_0 \cap C_1 = C_0$ so $(K_0, C_0) \subseteq (K_1, C_1)$ but the extension is not elementary since the formula $\exists u, v \in C (y=ux+v)$ (with parameters $x, y$) holds in $(K_1, C_1)$ but not in $(K_0,C_0)$. Note that this argument has been adapted from a standard proof of non-modularity of algebraically closed fields of transcendence degree at least $3$ (see, for example, \cite[Example 8.1.12]{Mar}).
\end{remark}

Now we are ready to prove Theorem \ref{thm:pairs-dichotomy} stating that any strongly minimal set definable in $\mathcal{K}_C$ is non-orthogonal to $C$.

\begin{proof}[Proof of Theorem \ref{thm:pairs-dichotomy}]
 {Let $S \subseteq K$ be a strongly minimal set defined over a finite tuple $\bar{a}$, and let $b \in S$ be generic over $\bar{a}$. Since $\Mr(b/\bar{a})=1$, by Lemma \ref{lem:pairs-MR-fin-alg} $b \in C(\bar{a})^{\alg}$. Proposition \ref{prop:pairs-type} asserts that $\tp(b/\bar{a})$ is determined by the formulae $\varphi_U(y)$. Hence, there is a cofinite subset $S'$ of $S$ which is defined by $\varphi_U$ for some $U$ since a conjunction of formulas of the form \eqref{eq-exist} is again of the same form. Then $S' \subseteq C(\bar{a})^{\alg} \subseteq \acl(C \cup \bar{a})$ and therefore $S' \not \perp C$. So $S \not \perp C$, for $S'$ is cofinite in $S$.}
\end{proof}

\begin{remark}
Let $S \subseteq K$ be strongly minimal defined by some formula $\varphi_U$.  As we pointed out above, $V:= W(K) \setminus U(K) \cap C^l$ is defined over $C$. So $V$ is a constructible set over $C$. Define an equivalence relation $E\subseteq V \times V$ by
$$ \bar{c}_1 E \bar{c}_2 \mbox{ iff } \forall y (p(\bar{a}, \bar{c}_1, y) = 0 \liff  p(\bar{a}, \bar{c}_2, y) = 0).$$

By stability $E$ is definable in the pure field structure of $C$. Moreover, there is a natural finite-to-one map from $S$ to $V/E$. By elimination of imaginaries in algebraically closed fields $V/E$ can be regarded as a constructible set in some Cartesian power $C^k$. The latter must have dimension $1$ since $S$ is strongly minimal. Thus, in the formula $\varphi_U$ we may assume that the constants live on a curve defined over $C$. This gives a characterisation of strongly minimal formulas. %Thus, we get a definable finite-to-one map from $S$ to a curve in $C^k$. Thus, strongly minimal sets in $\mathcal{K}_C$ are essentially algebraic curves over $C$.
\end{remark}

\subsection{Predimension for the differential equation of the $j$-function}\label{section-j}

Now we study the differential equation of the $j$-function. Subsections \ref{section-j} and \ref{section-EC} are preliminary. The reader is referred to \cite{Aslanyan-adequate-predim}  {and \cite{Aslanyan-Eterovic-Kirby-Diff-EC-j}} for details and proofs of the results.

 {In this section we let $\mathcal{K}:=(K;+,\cdot, D, 0, 1)$ be the countable saturated differentially closed field with field of constants $C$. This means, in particular, that $\td(C/\mathbb{Q}) = \aleph_0$.} Let $f_j(x,y)=0$ be the two-variable differential equation of the $j$-function (see Section \ref{jBackground}). We consider a binary predicate $E_j(x,y)$ which will be interpreted in a differential field as the set of solutions of the equation $f_j(x,y)=0$. This equation excludes the possibility of $x$ or $y$ being a constant. However, if we multiply $f_j(x,y)$ by a common denominator and make it a differential polynomial then $x$ and $y$ would be allowed to be constants as well. So we add $C^2$ to $E_j$, i.e. any pair of constants is in $E_j$. Further, define a relation $E_j^{\times}(x,y)$ by the formula $$E_j(x,y) \wedge x\notin C \wedge y\notin C.$$ %Note that $C$ is definable in terms of $E_j$ by the formula $E_j(0,x)$.

\begin{definition}
The $j$-reduct of $\mathcal{K}$ is the structure $\K:=(K;+, \cdot, E_j, 0, 1)$. 
\end{definition}

\begin{definition}
An $E_j$-\emph{field} is a substructure $\A$ of $\K$ with an algebraically closed underlying field containing $C$.   We let $\mathfrak{C}$ be the collection of all $E_j$-fields.
\end{definition}

\begin{remark}
We will identify $E_j$-fields with their domains and denote both the structure and the domain by $\A, \B, \X, \Y$, possibly with subscripts or superscripts. In particular, we identify $K$ with $\K$.

Note that the structure $(C;+, \cdot, E_j, 0, 1)$ is an $E_j$-field where $E_j$ is interpreted as $C^2$, and in fact it is the smallest $E_j$-field. Since its structure is that of a pure algebraically closed field, we will denote it by a standard letter $C$ (rather than $\C$).
\end{remark}

\begin{remark}
The formula $E_j(0,y)$ defines the field of constants in $\K$. Since by definition $E_j$-fields contain $C$, that formula defines $C$ in any $E_j$-field. This allows one to define the relation $E_j^{\times}$ in $E_j$-fields.
\end{remark}

\begin{definition}
If $\A$ is an $E_j$-field, then a tuple $\left(\bar{z},\bar{j}\right) \in \A^{2n}$ is called an $E_j$-\emph{point} if $\left(z_i, j_i\right) \in E_j(\A)$ for each $i=1,\ldots,n$. By abuse of notation, we let $E_j(\A)$ denote the set of all $E_j$-points in $\A^{2n}$ for any natural number $n$. 
\end{definition}

Now we describe the ``functional equations'' of $E_j$.

\begin{proposition}[{\cite[Lemmas 4.10, 4.11]{Aslanyan-adequate-predim}}]\label{funct-eq}
Let $\A$ be an $E_j$-field. 
\begin{itemize}%[label=\roman*]
    \item If $(z_i,j_i) \in E^{\times}_j(\A),~ i=1,2,$ and $\Phi_N(j_1, j_2)=0$ for some modular polynomial $\Phi_N$ then $z_2 = gz_1$ for some $g \in \SL_2(C)$.
    \item If $(z_1,j_1)\in E^{\times}_j(\A)$ and $(z_2,j_2)\in \A^2$ such that $\Phi_N(j_1,j_2)=0$ for some $\Phi_N$ and $z_2=gz_1$ for some $g\in \SL_2(C)$ then $(z_2,j_2)\in E^{\times}_j(\A)$.
\end{itemize}
\end{proposition}

The following is an immediate consequence of Theorem \ref{j-chapter-Ax-for-j}.

\begin{theorem}[Ax-Schanuel without derivatives]\label{j-Ax-without-der}
Let $\A$ be an $E_j$-field and let $(z_i,j_i)\in E_j^{\times}(\A), i =1,\ldots,n$. Then 
$$\td_CC(\bar{z},\bar{j}) \geq n+1,$$
unless for some $N$ and some $i\neq k$ we have $\Phi_N(j_i,j_k)=0$.
\end{theorem}

%It is easy to see that reducts of differential fields to the language $\mathfrak{L}_j$ are $E_j$-fields. %Note that $\mathbb{Q}^{\alg}$ is an $E_j$-field with $E_j(\mathbb{Q}^{\alg})=(\mathbb{Q}^{\alg})^4$ and it is the smallest $E_j$-field.

%Let $C$ be an algebraically closed field with $\td(C/\mathbb{Q})=\aleph_0$ and let $\mathfrak{C}$ be the collection of all $E_j$-fields $K$ with $C_K=C$.  Note that $C$ is an $E_j$-field with $E_j(C)=C^2$ and it is the smallest structure in $\mathfrak{C}$. From now on, by an $E_j$-field we understand a member of $\mathfrak{C}$. 

\begin{definition}
\begin{itemize}
    \item[]
    
    \item For a subset $S \subseteq K$ the $E_j$-\emph{closure} (or the $\mathfrak{C}$-\emph{closure}) of $S$ inside $\K$, denoted\footnote{In the differential setting $\langle S \rangle$ stands for the differential subfield generated by a set $S$. But from now on this notation will be used only for the $E_j$-closure.} $\langle S \rangle$, is the $E_j$-subfield of $\K$ generated by $S$, that is, the substructure of $\K$ with domain $C(S)^{\alg}$ and with the induced structure from $\K$. Similarly, for a tuple $\bar{s}\subseteq K$ the $E_j$-subfield generated by $\bar{s}$ is denoted by $\langle \bar{s} \rangle$.
    
    \item For $\A, \B \in \mathfrak{C}$ we write $\A\B$ for $\langle \A \cup \B \rangle$.
    
    \item A structure $\A \in \mathfrak{C}$ is \emph{finitely generated} if $\A=\langle S \rangle$ for some finite subset $S \subseteq \A$.
    
    \item The collection of all finitely generated structures from $\mathfrak{C}$ will be denoted by $\mathfrak{C}_{f.g.}$, and $\A\subseteq_{f.g.} \B$ means that $\A$ is a finitely generated substructure of $\B$.
\end{itemize}
\end{definition}

Note that $\mathfrak{C}_{f.g.}$ consists of those $E_j$-fields which have finite transcendence degree over $C$ (which, in fact, are not finitely generated as structures).

\begin{definition}
Let $\A, \B \in \mathfrak{C}_{f.g.}$ with $\A \seq \B$.
\begin{itemize}
    \item An $E_j$-\emph{basis of} $\B$ \emph{over} $\A$ is an $E_j$-point $\bar{b}=\left(\bar{z},\bar{j}\right)$ from $\B$ of maximal length satisfying the following conditions:
    \begin{itemize}
        \item $j_i$ and $j_k$ are modularly independent for all $i \neq k$,
        \item $\left(z_i, j_i\right) \notin \A^2$ for each $i$.
    \end{itemize}
    
    \item We let $\sigma(\B/\A)$ be the length of $\bar{j}$ in an $E_j$-basis of $\B$ over $\A$. Equivalently, an $E_j$-basis of $\B$ over $\A$ has length $2\sigma(\B/\A)$.
    
    \item  An $E_j$-\emph{basis} of $\B$ is an $E_j$-basis of $\B$ over $C$. We write $\sigma(\B)$ for $\sigma(\B/C)$.
    
    \item The \emph{relative predimension} of $\B$ over $\A$ is defined as  $$\delta(\B/\A)=\td(\B/\A)-\sigma(\B/\A).$$
    
    \item The \emph{predimension} of $\B$ is $$\delta(\B)=\delta(\B/C) = \td(\B/C) - \sigma(\B).$$
    
    %\item For finite sets (or tuples) $S, T \seq K$ we define $C,\mathcal{C},\mathtt{C},\mathsf{C}$
\end{itemize}
\end{definition}

Note that $\sigma$ is well defined by Proposition \ref{funct-eq}, and it is easy to see that for $\A \subseteq \B \in \mathfrak{C}_{f.g.}$ one has $\sigma(\B/\A)=\sigma(\B)-\sigma(\A)$. Hence, $\delta(\B/\A) = \delta(\B)-\delta(\A)$.

Note that the Ax-Schanuel inequality implies that $\sigma$ is finite for finitely generated $E_j$-fields. Moreover, for $\A, \B  \in \mathfrak{C}_{f.g.}$ the inequality $$\sigma(\A \B) \geq \sigma(\A)+\sigma(\B)-\sigma(\A \cap \B)$$ holds, where $\A\B$ denotes the $E_j$-field generated by $\A \cup \B$. Hence $\delta$ is \emph{submodular}, that is, for all $\A, \B \in \mathfrak{C}_{f.g.}$ we have  
$$\delta(\A\B)+ \delta(\A\cap \B) \leq \delta(\A) + \delta (\B).$$

In terms of the predimension the Ax-Schanuel inequality states exactly that $\delta(\A) \geq 0$ for all $\A \in \mathfrak{C}_{f.g.}$ with equality holding if and only if $\A=C$.

%One can also define the relative predimension $\delta(\B/\A)$ when $\A$ and $\B$ are not necessarily finitely generated, but $\td(B/A)$ is finite, but we will not need it here.

\begin{definition}
Let $\A, \B \in \mathfrak{C}$ with $\A \subseteq \B$. We say $\A$ is \emph{strong} (or \emph{self-sufficient}) in $\B$, and write $\A \leq \B$, if for all $\X \subseteq_{f.g.} \B$ we have $\delta(\X\cap \A) \leq \delta(\X)$. One also says $\B$ is a \emph{strong extension} of $\A$. An embedding $\A \hookrightarrow \B$ is strong if the image of $\A$ is strong in $\B$.
For a finite tuple $\bar{a} \subseteq \B$ we say $\bar{a}$ is strong in $\B$, and write $\bar{a}\leq \B$, if $\langle \bar{a} \rangle \leq \B$.
\end{definition}

If $\A \seq \B$ are finitely generated, then $\A \leq \B$ if and only if for any $E_j$-field $\X$ with $\A \seq \X \seq \B$ we have $\delta(\X/\A)\geq 0$.

\begin{definition}
For a set $S \subseteq K$ we define the \emph{self-sufficient closure} (or \emph{strong closure}) of $S$ in $\K$ by 
$$\lceil S \rceil := \bigcap_{\substack{\A \in \mathfrak{C}\\ \A \leq \K\\ S \subseteq \A \subseteq \K}}\A.$$
%Let $\lceil X \rceil:= \lceil X \rceil_{\mathbb{M}}$.
\end{definition}

By \cite[Lemma 2.12]{Aslanyan-adequate-predim}, $\lceil S \rceil \leq \K$. Note also that $\leq$ is transitive.

\begin{lemma}[{\cite[Lemma 2.14]{Aslanyan-adequate-predim}}]\label{sscl-min}
If $\X \subseteq_{f.g.} \K$ then
\begin{itemize}
\item  $\lceil \X\rceil$ is finitely generated, and

\item $\delta(\lceil \X \rceil)=\min\{\delta(\Y): \X\subseteq \Y \subseteq_{f.g.} \K\}.$
\end{itemize}

\end{lemma}

The predimension gives rise to a dimension function of a pregeometry in the following way. 
\begin{definition}
For $\X \subseteq_{f.g.} \K$ define
$$d(\X):= \min \{ \delta(\Y): \X \subseteq \Y \subseteq_{f.g.} \K \}=\delta(\lceil \X\rceil).$$
For a finite subset (or tuple) $S \subseteq \K$ set $d(S):=d(\langle S \rangle)$.

For $\X, \Y \subseteq_{f.g.} \K$ the relative dimension of $\Y$ over $\X$ is defined as $$d(\Y/\X):= d(\X\Y) - d(\X),$$ and for finite subsets $S, T \subseteq \K$ we define $$d(T/S):=d(\langle T \rangle/\langle S \rangle) = d(T\cup S) - d(S).$$
\end{definition}

Define an operator $\cl : \mathcal{P}(\K) \rightarrow \mathcal{P}(\K)$ (the latter is the power set of $\K$) as follows. For a finite set $S\subseteq \K$ let $$\cl(S)=\{ b \in \K: d(b/S)=0 \}$$ and for an infinite set $S\subseteq \K$ let
$$\cl(S) = \bigcup_{\substack{S_0 \subseteq S\\ |S_0|< \aleph_0}} \cl(S_0).$$
Then $(\K,\cl)$ is a pregeometry and $d$ is its dimension function. %We will drop $\K$ from the notation $d_{\K}$.

The class of $E_j$-fields with strong embeddings is a \emph{strong amalgamation class} (see \cite[Definition 2.18]{Aslanyan-adequate-predim}) which allows one to carry out an amalgamation-with-predimension construction (which is the uncollapsed version of a Hrushovski construction), and obtain a countable universal structure which is saturated with respect to strong embeddings (see \cite[Theorem 2.20]{Aslanyan-adequate-predim}). It follows from \cite[Theorem 4.40]{Aslanyan-adequate-predim} and \cite[Theorem 1.1]{Aslanyan-Eterovic-Kirby-Diff-EC-j} that this universal structure is isomorphic to $\K$. Thus, we obtain the following theorem.

\begin{theorem}
The $j$-reduct $\K$ of the countable saturated differentially closed field $\mathcal{K}$ enjoys the following properties.
\begin{itemize}
\item $\K$ is \emph{universal} with respect to strong embeddings, i.e. every $E_j$-field can be strongly embedded into $\K$.
\item $\K$ is \emph{saturated} with respect to strong embeddings, i.e. for every $\A, \B \in \mathfrak{C}_{f.g.}$ with strong embeddings $\A \hookrightarrow \K$ and $\A \hookrightarrow \B$ there is  a strong embedding of $\B$ into $\K$ over $\A$.
\end{itemize}

Furthermore, $\K$ is homogeneous with respect to strong substructures, that is, any isomorphism between finitely generated strong substructures of $\K$ can be extended to an automorphism of $\K$.
\end{theorem}

\subsection{Existential closedness}\label{section-EC}

A key property of $\K$ that will be used in the following subsections is \emph{Existential Closedness}. It states roughly that systems of equations in terms of $E_j$ always have solutions in $\K$ unless having a solution contradicts Ax-Schanuel. We will give a precise statement shortly.

\begin{definition}
Let $n$ be a positive integer, $k \leq n$ and $1\leq i_1 < \ldots < i_k \leq n$. Denote $\bar{i}=(i_1,\ldots,i_k)$ and define the projection map $\pr_{\bar{i}}:\K^{n} \rightarrow \K^{k}$ by
$$\pr_{\bar{i}}:(x_1,\ldots,x_n)\mapsto (x_{i_1},\ldots,x_{i_n}).$$
Further, define  $\ppr_{\bar{i}}:\K^{2n}\rightarrow \K^{2k}$ by
$$\ppr_{\bar{i}}:(\bar{x},\bar{y})\mapsto (\pr_{\bar{i}}\bar{x},\pr_{\bar{i}}\bar{y}).$$
\end{definition}

\begin{definition}
An irreducible algebraic variety $V \subseteq \K^{2n}$ is \emph{ {broad}} if for any $1\leq k\leq n$ and any $1\leq i_1 < \ldots < i_k \leq n$ we have $\dim \ppr_{\bar{i}} V \geq k$. We say $V$ is \emph{strongly  {broad}} if the strict inequality $\dim \ppr_{\bar{i}} V > k$ holds.
\end{definition}

\begin{lemma}[{\cite[Lemma 4.26]{Aslanyan-adequate-predim}}]\label{lemma-locus-basis-normal}
If $\A, \B \in \mathfrak{C}_{f.g.}$ with $\A\leq \B$ and $\bar{b} $ is an $E_j$-basis of $\B$ over $\A$ then $\Loc_A(\bar{b})$ is  {broad}.
\end{lemma}

\begin{theorem}[Existential Closedness, {\cite[Theorems 1.1 and 3.7]{Aslanyan-Eterovic-Kirby-Diff-EC-j}}]
\begin{itemize}
    \item[]
    
    \item[\emph{EC}] For each  {broad} variety $V \subseteq \K^{2n}$  the intersection $E_j(\K) \cap V(\K)$ is non-empty.
    
    \item[\emph{SEC}] For each  {broad} variety $V \subseteq \K^{2n}$ defined over a finite tuple $\bar{a}\subseteq \K$, the intersection $E_j(\K) \cap V(\K)$ contains a point Zariski generic in $V$ over $\bar{a}$.
\end{itemize}
\end{theorem}

EC and SEC stand for Existential Closedness and Strong Existential Closedness respectively. %Note that EC holds in all $j$-reducts of differentially closed fields, while SEC holds in saturated ones.

\subsection{Types in $\mathcal{\K}$}\label{section-j-types}

 {Recall that $\K=(K;+, \cdot, E_j, 0, 1)$ is the $j$-reduct of the countable saturated differentially closed field $\mathcal{K} = (K; +, \cdot, D, 0, 1)$, and we identify the set $K$ with the structure $\K$. Recall also that for a subset $S \seq \K$ the $E_j$-subfield of $\K$ with underlying field $C(S)^{\alg}$ is denoted by $\langle S \rangle$.}% and $\delta(S)$ and $\sigma(S)$ stand for $\delta(\langle S \rangle)$ and $\sigma(\langle S \rangle)$ respectively.} %Further, as in the previous subsections, $\A, \B, \B'$ will always denote $E_j$-fields with underlying fields $A, B, B'$ respectively.}

\begin{lemma}\label{lemma-closure-basis-generated}
Let $\bar{a}\subseteq \K$. If $(\bar{u},\bar{v})$ is an $E_j$-basis of $\lceil \bar{a} \rceil$ then the latter is generated by $\bar{a}, \bar{u}, \bar{v}$,  {that is,  $\lceil \bar{a} \rceil = \langle \bar{a},\bar{u},\bar{v}\rangle.$ }
\end{lemma}
\begin{proof}
Let $\A := \langle \bar{a},\bar{u},\bar{v}\rangle = C(\bar{a}, \bar{u}, \bar{v})^{\alg} \subseteq  \lceil \bar{a} \rceil$. Then $\td(\lceil \bar{a} \rceil/C)\geq \td(\A/C)$ and $\sigma(\lceil \bar{a} \rceil) = |\bar{v}| = \sigma(\A)$. Hence $\delta(\A) \leq \delta(\lceil \bar{a} \rceil)$ and, by Lemma \ref{sscl-min}, $\delta(\A) = \delta(\lceil \bar{a} \rceil)$. Therefore $\td(\lceil \bar{a} \rceil/C) = \td(\A/C)$ and $\A = \lceil \bar{a} \rceil.$
\end{proof}

Let $\bar{a}=(a_1,\ldots,a_m)\in \K^m$ be a tuple with $d(\bar{a})=k$ and $b\in \K$ with $d(b/\bar{a})=0$, i.e. $b \in \cl (\bar{a})$. This means that $d(\bar{a}b)=k$. Pick an $E_j$-basis $(\bar{z},\bar{j})$ of $\B:=\lceil \bar{a}, b \rceil$. 
By Lemma \ref{lemma-closure-basis-generated} $\B=C(\bar{a},b, \bar{z},\bar{j})^{\alg}$. We claim that $b \in C(\bar{a}, \bar{z},\bar{j})^{\alg}$. Indeed, if it is not true then $$k=d(\bar{a}b)= \delta (\B) > \delta ( C(\bar{a}, \bar{z},\bar{j})^{\alg}) \geq d(\bar{a}) = k .$$

Thus, $\B=C(\bar{a}, \bar{z},\bar{j})^{\alg}$. Let $p(\bar{a},\bar{c},\bar{z},\bar{j},b)=0$ for some irreducible polynomial $p$ and $\bar{c}\in C^s$. Denote $l:=|\bar{j}| = \sigma(\B)$ and $V:= \Loc_{C(\bar{a})}(\bar{z},\bar{j}) \subseteq \K^{2l}$ and assume (by extending $\bar{c}$ if necessary) it is defined over $\bar{c},\bar{a}$. In order to stress that $V$ is defined over $\bar{a},\bar{c}$, we denote it by $V_{\bar{c}}$ (we do not include $\bar{a}$ in this notation, for it is fixed and will not vary). This way we get a parametric family of varieties as $\bar{c}$ varies. Notice that 
\begin{equation}\label{eq-1}
 \dim V_{\bar{c}} = \td (\B/C) - \td(\bar{a}/C) = \delta(\B) + \sigma(\B) - \td(\bar{a}/C) = k+l - \td(\bar{a}/C).   %= \td(B/A) = \delta(B/A) - \sigma(B/A). %
\end{equation}

Also, let $W:=\Loc_{\mathbb{Q}(\bar{a})}(\bar{c})$. For each proper Zariski closed subvariety $U\subsetneq W$, defined over $\mathbb{Q}(\bar{a})$, and for each positive integer $N$ consider the formulae
\begin{gather*}\label{equation-type-E-formula}
    \xi_{U,N}(\bar{e},\bar{u},\bar{v}):= \left( \bar{e}\in C^{|\bar{c}|} \cap (W\setminus U) \wedge   (\bar{u},\bar{v})\in V_{\bar{e}} \cap E^{\times}_j \wedge \bigwedge_{n=1}^N \bigwedge_{i\neq r} \Phi_n(v_i,v_r)\neq 0 \right),\\
    \psi_{U, N} (\bar{e},\bar{u},\bar{v},y) :=  \xi_{U,N}(\bar{e},\bar{u},\bar{v})  \wedge p(\bar{a},\bar{e},\bar{u},\bar{v},y) = 0 ,\\
     \varphi_{U,N}(y):= \exists \bar{e} , \bar{u}, \bar{v}~ \psi_{U,N}(\bar{e},\bar{u},\bar{v},y).
\end{gather*}

Observe that $\varphi_{U,N}$ is defined over $\bar{a}$ and $\varphi_{U,N}(b)$ holds in $\K$. Obviously these formulas depend on $\bar{a}$ and $b$, but it is not explicit in their notation, for it will always be clear what $\bar{a}$ and $b$ are. Thus, in this and the following subsections we will use the notation introduced above for given $\bar{a}$ and $b$.

\begin{proposition}\label{prop-types-axiom}
The formulae $\varphi_{U,N}$ axiomatise the type of $b$ over $\bar{a}$ in $\K$, that is, if for some $b'\in \K$ the formula $\varphi_{U,N}(b')$ holds in $\K$ for each $U\subsetneq W$ and each $N>0$ then $\tp(b/\bar{a}) = \tp(b'/\bar{a})$.
\end{proposition}
\begin{proof}
Consider the type $q(\bar{e},\bar{u},\bar{v})$ over $\bar{a},b'$ consisting of all formulae $\psi_{U, N} (\bar{e},\bar{u},\bar{v},b')$ for all $U\subsetneq V$ and $N>0$. Then $q$ is finitely satisfiable and hence there is a realisation of $q$ in $\K$ which we denote by $(\bar{c}',\bar{z}',\bar{j}')$.

Observe that $\bar{c}'$ is generic in $W$ over $\bar{a}$. Hence $\bar{c}$ and $\bar{c}'$ have the same algebraic type over $\bar{a}$. %there is an automorphism of the field $C(\bar{a})^{\alg}$ which fixes $\bar{a}$ pointwise, fixes $C$ setwise and sends $\bar{c}$ to $\bar{c}'$. 
In particular, $\dim V_{\bar{c}} = \dim V_{\bar{c}'}$. Further, $j'_1,\ldots,j'_l$ are pairwise modularly independent. So if $\B':=C(\bar{a},\bar{z}',\bar{j}')^{\alg}$ then $\sigma(\B')\geq l$. On the other hand,
\begin{equation}\label{eq-inequality}
   \td(\B'/C) = \td(\bar{a}/C) + \td(\bar{z}',\bar{j}'/C(\bar{a})) \leq \td(\bar{a}/C) + \dim V_{\bar{c}'} = k+l,
\end{equation}
where the last equality follows from \eqref{eq-1}. Therefore $\delta(\B') \leq (k+l) -l = k$. However, $\B'$ contains $\bar{a}$ and since $d(\bar{a}) = k$, $\delta(\B')$ cannot be smaller than $k$. Thus, $\delta(\B')=k$ and $\sigma(\B')=l$ and the inequality in \eqref{eq-inequality} must be an equality, i.e. $\td(\bar{z}',\bar{j}'/C(\bar{a})) =  \dim V_{\bar{c}'}$. This means that $(\bar{z}',\bar{j}')$ is generic in $V_{\bar{c}'} $ over $C(\bar{a})$. Therefore, there is a field isomorphism $\pi:\B\rightarrow \B'$ which fixes $\bar{a}$ pointwise, fixes $C$ setwise, sends $\bar{c}$ to $\bar{c}'$ and sends $(\bar{z},\bar{j})$ to $(\bar{z}',\bar{j}')$. %the map $\pi$ can be extended to a field isomorphism $\pi:B\rightarrow B'$ which sends $(\bar{z},\bar{j})$ to $(\bar{z}',\bar{j}')$. 
Since $\sigma(\B')=l$, the tuple $(\bar{z}',\bar{j}')$ is an $E_j$-basis of $\B'$ and $\pi$ is an isomorphism of $\B$ and $\B'$ as $E_j$-fields.

Finally, as $p(\bar{a},\bar{c},\bar{z},\bar{j},b) = 0$ and $p(\bar{a},\bar{c}',\bar{z}',\bar{j}',b') = 0$, we could have chosen $\pi$ so that $\pi(b) = \pi(b')$. Now both $\B$ and $\B'$ are strong in $\K$ and the latter is homogeneous with respect to strong substructures, hence $\pi$ can be extended to an automorphism of $\K$. This shows that $b$ and $b'$ have the same type in $\K$ over $\bar{a}$.
\end{proof}

\begin{remark}
In general, all types in $\K$ are determined by formulas of the form $\varphi_{U,N}$ and their negations. In particular, every formula is equivalent to a Boolean combination of existential formulas in $\K$ and hence its theory is near model complete. In \cite[\S 4]{Aslanyan-adequate-predim} we gave an axiomatisation of $\Th(K)$. %Moreover, ``small'' sets, i.e. sets of finite Morley rank, are existentially definable while ``large'' sets, i.e. sets of Morley rank $\omega$, are universally definable.
\end{remark}

\begin{theorem}\label{theorem-j-acl}
If $\bar{a}\leq \K$ then $\acl(C(\bar{a})) = C(\bar{a})^{\alg}.$
%$$\acl (A) = \mathbb{Q}(A)^{\alg}.$$
\end{theorem}
\begin{proof}
It suffices to prove that for $\bar{a}\leq \K$ we have $\acl(\bar{a})\subseteq C(\bar{a})^{\alg}$. Assume $b\in \acl(\bar{a})$. Then $d(b/\bar{a})=0$ and $\tp(b/\bar{a})$ is determined by existential formulas $\varphi_{U,N}(y)$. Since $b\in \acl(\bar{a})$, some formula $\varphi_{U,N}(y) \in \tp(b/\bar{a})$ has finitely many realisations in $\K$. %Denote those realisations by $b_1,\ldots, b_t$.

We use the notation introduced above. The point $(\bar{z},\bar{j})\in V_{\bar{c}}$ is generic over $\bar{a},\bar{c}$. Observe that $(\bar{z},\bar{j})$ must contain an $E_j$-basis of $\A:=\langle \bar{a} \rangle$. Denote it by $(\bar{z}_{\bar{a}}, \bar{j}_{\bar{a}})$ and let $(\bar{z}_0, \bar{j}_0) := (\bar{z},\bar{j}) \setminus (\bar{z}_{\bar{a}}, \bar{j}_{\bar{a}})$, i.e. $(\bar{z}_0,\bar{j}_0)$ consists of all coordinates $(z_i,j_i)$ of $(\bar{z},\bar{j})$ for which $(z_i,j_i)\notin \A^2$. In other words, $(\bar{z}_0,\bar{j}_0)$ is an $E_j$-basis of $\B:=\lceil \A(b) \rceil$ over $\A$. Let $W$ be an irreducible component over $L:=\mathbb{Q}(\bar{a},\bar{c},\bar{z}_{\bar{a}}, \bar{j}_{\bar{a}})^{\alg}$ of the fibre of $V_{\bar{c}}$ above $(\bar{z}_{\bar{a}}, \bar{j}_{\bar{a}})$ containing $(\bar{z}_0, \bar{j}_0)$. Then it is defined over $L$ and $(\bar{z}_0, \bar{j}_0)$ is generic in $W$ over $L$. 

Since $\varphi_{U,N}(b)$ holds in $\K$, in particular we have $p(\bar{a},\bar{c}, \bar{z}, \bar{j}, b)=0$. Assume
\begin{equation*}
    p(\bar{a},\bar{c}, \bar{z}, \bar{j}, Y)= Y^n + s_{n-1}( \bar{z}_0, \bar{j}_0) Y^{n-1}+\cdots+ s_{0}(\bar{z}_0, \bar{j}_0)
\end{equation*}
where each $s_i(\bar{X}_1,\bar{X}_2)$ is a rational function over $L$. If $s_i(\bar{z}_0, \bar{j}_0) \in L$ for all $i$, then $b\in L\subseteq \A$ and we are done. Otherwise assume without loss of generality that $s_0(\bar{z}_0, \bar{j}_0)\notin L$.

Since $\A\leq \B$, by Lemma \ref{lemma-locus-basis-normal}, $W$ is  {broad}. By SEC there is a point $(\bar{z}_1,\bar{j}_1)\in W(\K)\cap E_j(\K)$ generic over $L(\bar{z}_0,\bar{j}_0)$. If $s_0(\bar{z}_1, \bar{j}_1) = s_0(\bar{z}_0, \bar{j}_0)$ then the function $s_0(\bar{X}_1,\bar{X}_2)$ is constant on $W$. On the other hand, $W$ is defined over $L$, so the constant value of $s_0(\bar{X}_1,\bar{X}_2)$ must belong to $L$. This is a contradiction, hence  $s_0(\bar{z}_1, \bar{j}_1) \neq s_0(\bar{z}_0, \bar{j}_0)$. Now pick a generic point $(\bar{z}_2,\bar{j}_2)$ in $W(\K)\cap E_j(\K)$ over $L(\bar{z}_0,\bar{j}_0,\bar{z}_1,\bar{j}_1)$. By the above argument the elements $s_0(\bar{z}_0, \bar{j}_0), s_0(\bar{z}_1, \bar{j}_1), s_0(\bar{z}_2, \bar{j}_2)$ are pairwise distinct. Repeating this procedure countably many times we construct a sequence $(\bar{z}_i, \bar{j}_i), i=0,1, 2, \ldots$, such that for each $i$ $$\K\models \xi_{U,N}(\bar{c},\bar{z}_{\bar{a}},\bar{z}_i,\bar{j}_{\bar{a}},\bar{j}_i)$$ and $s_0(\bar{z}_i, \bar{j}_i), i=0,1,2,\ldots$, are pairwise distinct. This shows that the formula $\varphi_{U,N}(y)$ has infinitely many realisations, for there are only finitely many monic polynomials of a given degree the roots of which belong to a finite set of elements. This is a contradiction.
\end{proof}

\begin{corollary}
For any $\bar{a}\subseteq \K$ we have $\acl(\bar{a})\subseteq \lceil \bar{a} \rceil$.
\end{corollary}

\subsection{Classification of strongly minimal sets in $\K$}

 {We are now ready to prove Theorem \ref{thm:intro-K-dichotomy} which we restate for convenience. Recall that we still work in the $j$-reduct $\K=(K;+, \cdot, E_j, 0, 1)$ of the countable saturated differentially closed field $\mathcal{K} = (K; +, \cdot, D, 0, 1)$.}

\begin{customthm}{\ref{thm:intro-K-dichotomy}}
Let $S\subseteq \K$ be a strongly minimal set. Then either $S$ is geometrically trivial or $S \not \perp C$.
\end{customthm}

\begin{proof}
Assume $S$ is defined over $\bar{a}$ and $S \perp C$. Let $\A:=\langle \bar{a} \rangle$ be the $E_j$-closure of $\bar{a}$.  {First, we extend $\A$ to a possibly larger $E_j$-field and work over it which simplifies some technical arguments.} Pick $b\in S\setminus \acl(\A)$ (if such an element does not exist then $S\not \perp C$). Then clearly $d(b/\A)=0$. Let $\B':=\lceil \A b \rceil$ and let $z_b\in \K$ be such that $E_j(z_b, b)$ holds. Now if $\B:= \langle B'(z_b)\rangle$ then $\delta(\B) = \delta(\B') = d(\A)$ as $d(b/\A)=0$. Hence $\B\leq \K$. Choose a maximal $E_j$-field $\A'$ with $\A\subseteq \A' \subseteq \B$ such that $b\notin \acl(\A')$. Since strong minimality and geometric triviality of $S$ do not depend on the choice of the set of parameters over which $S$ is defined (see Theorem \ref{geometric-triviality}), we may extend $\A$ and assume $\A'=\A$. This means that if $e\in \B\setminus \A$ then $b\in \acl(\A e)$. In particular, $\acl(\A)=\A$.

 {We will show that any pairwise $\acl$-independent elements from $S$ over $\A$ are independent over $\A$. It is equivalent to geometric triviality by Lemma \ref{lem:append-geom-triv}.} 

Let $(\bar{z},\bar{j})\in \B^{2l}$ be an $E_j$-basis of $\B$ with $j_l=b$. By extending $\bar{a}$ we may assume that $V:=\Loc_{\A}(\bar{z},\bar{j})\subseteq \K^{2l}$ is defined over $\bar{a}$. By the proof of Proposition \ref{prop-types-axiom}, $\tp(b/\A)$ is determined by the formulae 
$$\chi_N(y) := \exists \bar{u}, \bar{v} \left( (\bar{u},\bar{v})\in V \cap E^{\times}_j  \wedge \bigwedge_{n=1}^N \bigwedge_{i\neq r} \Phi_n(v_i,v_r)\neq 0 \wedge y=v_l \right).$$

Now pick pairwise $\acl$-independent elements $b_1,\ldots, b_t \in S\setminus \A$. We need to prove that $b_t \notin \acl(\A b_1\ldots b_{t-1})$. Since $S$ is strongly minimal, $\tp(b_i/\A)=\tp(b/\A)$ for all $i$. By saturatedness of $\K$ for each $i$ there is $(\bar{z}^i, \bar{j}^i)\in V\cap E_j^{\times}$ such that $\bar{j}^i$ is pairwise modularly independent\footnote{This means that the coordinates of $\bar{j}^i$ are pairwise modularly independent.} and $j^i_l=b_i$. Let $\B_i:=\langle \A(\bar{z}^i,\bar{j}^i)\rangle$.

It is clear that $\dim V = \td(\bar{z},\bar{j}/\A) = \td(\B/A) = \delta(\B/\A) + \sigma(\B/\A)$. Therefore
\begin{gather*}
    \delta(\B_i) = \td(\B_i/C) - \sigma(\B_i) \leq \dim V + \td(\A/C) - l = \\ \delta(\B/\A) + \sigma(\B/\A) + \delta(\A) +\sigma(\A) - l = \delta(\B) = d(\A),
\end{gather*}
and so $\B_i \leq \K$ and $(\bar{z}^i, \bar{j}^i)$ is an $E_j$-basis of $\B_i$. We can conclude now that $\lceil \A b_i \rceil \subseteq \B_i$, hence $\acl(\A b_i)\subseteq \B_i$.  Moreover, as in the proof of Proposition \ref{prop-types-axiom} there is an automorphism of $\K$ over $\A$ that maps $\B$ onto $\B_i$ (and maps $(\bar{z},\bar{j})$ to $(\bar{z}^i,\bar{j}^i)$). In particular, for every $e\in \B_i\setminus \A$ we have $b_i \in \acl(\A e)$.

We claim that $j^i_r$ and $j^m_k$ are modularly independent unless $(i,r)=(m,k)$ or $j^i_r, j^m_k \in \A$. Assume for contradiction that for some $i\neq m$ the elements $j^i_r$ and $j^m_k$ are modularly dependent and $j^i_r\notin \A$. Then $b_i \in \acl(\A j^i_r) = \acl(\A j^m_k) \subseteq \B_m$. Hence $b_m \in \acl(\A b_i)$ which is a contradiction, for we assumed $b_i$'s are pairwise $\acl$-independent over $\A$. This shows, in particular, that (when $t\geq 2$) $\A\leq \K$ as otherwise we would have $b\in \lceil \A \rceil$ and $S \subseteq \acl(\A b)$ in which case $S \not \perp C$.

Now let $\tilde{\B}_k:= \langle \B_1 \cup \ldots \cup \B_k \rangle$ be the $E_j$-subfield of $\K$ generated by $\B_1 \cup \ldots \cup \B_k$ where $k\leq t$. The above argument shows that $$\sigma(\tilde{\B}_k/\A) = \sum_{i=1}^k \sigma(\B_i/\A) = k\cdot \sigma(\B/\A) . $$

By submodularity of $\delta$ we have
$$\delta(\tilde{\B}_k)\leq \delta(\tilde{\B}_{k-1})+\delta(\B_k) - \delta(\tilde{\B}_{k-1}\cap \B_k)$$ for each $k$. Since $\delta(\tilde{\B}_{k-1}\cap \B_k)\geq d(\A)$, we can show by induction that $\delta(\tilde{\B}_k)=d(\A)$ and $\tilde{\B}_k\leq \K$. Thus, 
\begin{equation}\label{eq-td-delta}
  \td(\tilde{\B}_k/C) = \delta (\tilde{\B}_k) + \sigma(\tilde{\B}_k) = d(\A)  + \sigma(\tilde{\B}_k).  
\end{equation}

On the other hand, using submodularity of $\td$ and $-\sigma$ we get by induction
\begin{gather*}
    \td(\tilde{\B}_k/C) \leq \td(\tilde{\B}_{k-1}/C) + \td(\B_k/C) - \td((\tilde{\B}_{k-1} \cap \B_k)/C) = \\
    d(\A) + \sigma(\tilde{\B}_{k-1}) + d(\A) + \sigma(\B_k) - \delta(\tilde{\B}_{k-1} \cap \B_k) - \sigma(\tilde{\B}_{k-1} \cap \B_k) \\ \leq d(\A) + \sigma(\tilde{\B}_k),
\end{gather*}
where $\delta(\tilde{\B}_{k-1} \cap \B_k) \geq d(\A)$ for $\A\subseteq \tilde{\B}_{k-1} \cap \B_k$. In fact we must have equalities everywhere in the above inequality due to \eqref{eq-td-delta}. In particular, 
$$\sigma((\tilde{\B}_{k-1}\cap \B_k)/\A) = \sigma(\tilde{\B}_{k-1}/\A)+ \sigma (\B_k/\A) - \sigma (\tilde{\B}_k/\A) = 0.$$

So $$\td((\tilde{\B}_{k-1}\cap \B_k)/\A) = \delta ((\tilde{\B}_{k-1}\cap \B_k)/\A) + \sigma((\tilde{\B}_{k-1}\cap \B_k)/\A) = 0.$$ This implies that $\tilde{\B}_{k-1}\cap \B_k=\A$. In particular, $b_t \notin \tilde{\B}_{t-1}$. On the other hand, $\acl(\A b_1\ldots b_{t-1}) \subseteq \lceil \A b_1\ldots b_{t-1} \rceil \subseteq \tilde{\B}_{t-1}$. Thus, $b_t \notin \acl(\A b_1\ldots b_{t-1})$ as required.
\end{proof}

We can also prove that certain sets are strongly minimal. Let $\A:=C(\bar{a})^{\alg}\leq \K$. Assume $V\subseteq \K^2$ is an algebraic curve defined over $\A$, i.e. $\dim V=1$. Consider the formula
$$\chi(y):= \exists \bar{u}, \bar{v} \left( (\bar{u},\bar{v})\in V \cap E^{\times}_j \wedge  p(\bar{a},\bar{u},\bar{v},y) = 0 \right), $$ where $p$ is some irreducible algebraic polynomial.

%Note that $\chi_N (y) = \exists \bar{u}, \bar{v} \psi_{U,N}(\bar{c}, \bar{u}, \bar{v}, y)$ is defined over $\bar{a}, \bar{c}$.

\begin{proposition}\label{prop-SM-classification}
If $S:=\chi(\K)$ is infinite then $S$ is strongly minimal.
\end{proposition}

\begin{proof}

We need to show that over any set of parameters all non-algebraic elements in $S$ realise the same type. By stability we may choose all extra parameters from the set $S$ itself. Assume $e, e', b_1, \ldots, b_t \in S$ with $e, e' \notin \A( \bar{b})^{\alg}$. We will show that $\tp(e/\A( \bar{b}))=\tp(e'/\A(\bar{b}))$.

Choose existential witnesses $({z}, {j}), ({z}', {j}'), ({z}_i, {j}_i)\in V(\K)\cap E_j^{\times}(\K)$ for $\chi(e),~ \chi(e')$ and $\chi(b_i)$ respectively. Since $e \notin \A(\bar{b})^{\alg}$ and $p(\bar{a},{z},{j},e)=0$ and $\dim V=1$, the point $({z}, {j})$ is generic in $V$ over $\A(\bar{b})$. Similarly $({z}', {j}')$ is generic in $V$. So $(z,j)$ and $(z',j')$ have the same algebraic type over $\A(\bar{b})$. On the other hand, $\delta(\bar{b}/\A)\leq 0$, therefore $\delta(\bar{b}/\A) = 0$. Thus $\delta(e/\A\bar{b}) = \delta(e'/\A\bar{b})$ and $(z, j)$ and $(z',j')$ form $E_j$-bases of $\A(\bar{b}, e)^{\alg}$ and $\A(\bar{b}, e')^{\alg}$ over $\A(\bar{b})^{\alg}$ respectively. Hence, as in the proof of Proposition \ref{prop-types-axiom},  $e$ can be mapped to $e'$ by an automorphism of $\K$ over $\A(\bar{b})$.
\end{proof}

\begin{remark}
When $\A$ is not strong in $\K$ we may actually work over $\lceil \A \rceil$ since strong minimality of a set does not depend on the choice of the set of parameters over which it is defined. Hence the assumption $\A\leq \K$ does not restrict generality.
\end{remark}

\appendix
\section{On strong minimality}\label{A}

In this appendix we give some preliminaries on strongly
minimal sets. %that are used throughout the paper. %We also prove that geometric triviality of a strongly minimal set does not depend on the set of parameters over which our set is defined. It is of course a well known classical result, but we sketch a proof here since we use the result in the proof of Theorem \ref{cor}. 
For a detailed account of strongly minimal sets and geometric stability theory in general we refer the reader to \cite{Pillay-geometric}. 

Algebraic closure defines a pregeometry on a strongly minimal set. More precisely, if $X$ is a strongly minimal set in a structure $\mathcal{M}$ defined over $A \subseteq M$ then the operator
$$\cl_A : Y \mapsto \acl(AY) \cap X, \text{ for } Y \subseteq X,$$
is a pregeometry.  {This depends on the set of parameters $A$ and, in particular, we get a different pregeometry if we extend $A$. Nevertheless, many important properties of the pregeometry are independent of the choice of the parameter set. Theorem \ref{geometric-triviality} is an example of that.} %When the set of parameters $A$ does not play any role we write $\cl$ instead of $\cl_A$.}

\begin{definition}\label{geometry} 
Let $\mathcal{M}$ be a structure (with domain $M$) and $X \subseteq M$ be a strongly minimal set defined over a finite set $A \subseteq M$.
\begin{itemize}
\item We say $X$ is \emph{geometrically trivial} (over $A$) if whenever $Y \subseteq X$ and $z \in \cl_A(Y)$, we have $z \in \cl_A(y)$ for some $y \in Y$. In other words, geometric triviality means that the closure of a set is equal to the union of closures of singletons.
\item $X$ is called \emph{strictly disintegrated} (over $A$) if any distinct elements $x_1, \ldots, x_n \in X$ are independent (over $A$).
\item $X$ is called $\aleph_0$-\emph{categorical} (over $A$) if it realises only finitely many $1$-types over $AY$ for any finite $Y \subseteq X$. This is equivalent to saying that  $\cl_A(Y)$ is finite for any finite $Y \subseteq X$.
\end{itemize}
\end{definition}

Note that strict disintegratedness implies $\aleph_0$-categoricity and geometric triviality.

\begin{lemma}\label{lem:append-geom-triv}
 {Let $X$ be a strongly minimal set. Then $X$ is geometrically trivial if and only if any finitely many pairwise independent elements of $X$ are independent.}
\end{lemma}
\begin{proof}
 {Here we work over a parameter set $A$ and independent means $\cl_A$-independent. Since $A$ is fixed, we will write $\cl$ instead of $\cl_A$.} 

 {Assume $X$ is geometrically trivial and let $x_1,\ldots,x_n \in X$ be dependent. Then we may assume without loss of generality that $x_n \in \cl(x_1,\ldots,x_{n-1})$. Hence by triviality $x_n\in \cl(x_i)$ for some $i\in \{1,\ldots,n-1\}$. So $x_1,\ldots,x_n$ are not pairwise independent.}

 {Conversely, assume any pairwise independent elements from $X$ are independent and let $z\in \cl(Y)$ for some $Y \seq X$. We may assume $Y$ is finite and independent. Then $Y \cup \{ z\}$ is a dependent set, so it is not pairwise independent. Therefore $z \in \cl(y)$ for some $y\in Y$.}
\end{proof}

\begin{theorem}[{cf. \cite[Lemma 2.20]{Nagloo-Pillay0}}]\label{geometric-triviality}
Let $\mathcal{M}$ be a model of an $\omega$-stable theory and $X \subseteq M$ be a strongly minimal set which is definable over $A \seq M$ and also over $B\seq M$. Then $X$ is geometrically trivial over $A$ if and only if it is geometrically trivial over $B$.
\end{theorem}
\begin{proof}
We may assume $A$  and $B$ are finite and $A\seq B$. Moreover, by expanding the language with constant symbols for elements of $A$ we may assume further that $X$ is $\emptyset$-definable, i.e. $A  = \emptyset$. Let $\bar{b}:=(b_1,\ldots,b_m)$ be an enumeration of $B$. Let $z \in \acl(BY)$ for some finite $Y \subseteq X$. By stability $\tp(\bar{b}/X)$ is definable over a finite $C \subseteq X$ and we may assume that $C \subseteq \acl(B) \cap X$.  Therefore $z \in \acl(CY)$. By geometric triviality of $X$ over $\emptyset$ we have $z \in \acl(c)$ for some $c \in C$ or $ z \in \acl(y)$ for some $  y\in Y$. This shows geometric triviality of $X$ over $B$.

Conversely, assume $X$ is not geometrically trivial over $\emptyset$. Then there are elements $y_1,\ldots,y_n \in X$ which are dependent but pairwise independent. Let $\bar{z}\in X^n$ be a realisation of a non-forking extension of $\tp(\bar{y})$ to $B$. Then clearly $\bar{z}$ is $\cl_B$-dependent but pairwise $\cl_B$-independent which means $X$ is not geometrically trivial over $B$. 
\end{proof}

As we saw in the proof, all definable subsets of $X^n$ over $B$ are definable over $\acl(B) \cap X$ (which means that $X$ is \emph{stably embedded} into $M$). The same argument shows that $\aleph_0$-categoricity does not depend on parameters (see \cite[Lemma 2.20]{Nagloo-Pillay0}). Of course this is not true for strict disintegratedness but a weaker property is preserved. Namely, if $X$ is strictly disintegrated over $A$ then any distinct non-algebraic elements over $B$ are independent over $B$.

\iffalse
\begin{definition}
Two definable sets $X$ and $Y$ are called \emph{orthogonal}, written $X \perp Y$, if any two elements $x\in X$ and $y\in Y$ are (forking) independent over any set of parameters over which $X$ and $Y$ are defined.
\end{definition}

The following gives a simpler characterisation of orthogonality for strongly minimal sets.

\begin{lemma}
Two strongly minimal sets $X$ and $Y$ are non-orthogonal if and only if for some finite parameter set $A$ we have $Y\subseteq \acl(A\cup X)$.
\end{lemma}
\fi

\begin{definition}
Let $X$ and $Y$ be strongly minimal sets definable in a structure $\mathcal{M}$. Then $X$ and $Y$ are called \emph{non-orthogonal}, written $X \not \perp Y$, if for some finite set $A\seq M$ we have $Y\subseteq \acl(A\cup X)$.
\end{definition}

Non-orthogonality means that the given sets are ``similar''. It is an equivalence relation for strongly minimal sets.

\subsection*{Acknowledgements} The first part of this work was done while I was a DPhil student at the University of Oxford under the supervision of Boris Zilber and Jonathan Pila. I would like to thank them for valuable suggestions and comments. I am also grateful to the referee for a thorough reading of the paper and numerous helpful comments.

%\addcontentsline {toc} {section} {Bibliography}
\bibliographystyle {alpha}
\bibliography {ref}

\end{document}